\documentclass[12pt]{article}

\usepackage{amsmath}
\usepackage{amssymb}
\usepackage{graphicx}

\newtheorem{thr}{Theorem}[section]
\newtheorem{lm}{Lemma}[section]
\newtheorem{assum}{Assumption}
\newtheorem{prop}{Proposition}[section]

\usepackage{verbatim}
\usepackage{color}

\newcommand{\bea}{\begin{eqnarray}}
\newcommand{\eea}{\end{eqnarray}}
\newcommand{\ri}{\rightarrow}

\newcommand{\E}{\mathsf{E}}
\newcommand{\Var}{\mathsf{Var}}
\newcommand{\Cov}{\mathsf{Cov}}
\renewcommand{\P}{\mathsf{P}}
\newcommand{\R}{{\textbf{R}}}
\newcommand{\Cc}{{\textbf{C}}}
\newcommand{\Z}{{\textbf{Z}}}

\newcommand{\by}{{\bf y}}
\newcommand{\bx}{{\bf x}}

\newcommand{\SSS}{{\cal S}}

\newcommand{\bean}{\begin{eqnarray*}}
\newcommand{\eean}{\end{eqnarray*}}
\newcommand{\nn}{\nu}

\newcommand{\xbi}{\bar{\xi}_i}
\newcommand{\xbj}{\bar{\xi}_j}

\newcommand{\toP}{\stackrel{{\P_{m,\beta}}}{\longrightarrow}}
\newcommand{\toPP}{\stackrel{{\P^{(0)}_{m}}}{\longrightarrow}}

\begin{document}

\title{Asymptotic normality of maximum likelihood estimator 
for
cooperative sequential adsorption\thanks{This research was supported by the 
Royal Society of London, Travel for Collaboration Grant,
 award ref. code RC-MA1038}}

\author{Mathew D. Penrose\footnote{Postal address: Department of Mathematical Sciences,
University of Bath, BA2 7AY, UK.
 Email address: masmdp@bath.ac.uk}\\
{\small University of Bath}
\and 
Vadim Shcherbakov\footnote{Postal address: Laboratory of Large Random Systems,
 Faculty of Mechanics and Mathematics, Moscow
State University, Glavnoe Zdanie, Leninskie Gory, Moscow, 119991, Russia. 
Email address: v.shcherbakov@mech.math.msu.su}\\
{\small  Moscow State University}}

\maketitle

\vspace*{-5ex}

\begin{abstract}

We have shown in previous work that statistical inference 
for cooperative sequential adsorption  model 
can be based on maximum likelihood estimation.
In this paper we continue this research 
and establish  asymptotic normality of the maximum likelihood 
 estimator in thermodynamic limit.
We also perform and discuss some numerical simulations of the model.


\end{abstract}

\section{Introduction}
\label{intro}

 This paper continues  the research started in 
\cite{MPVS}, where  properties  of maximum 
likelihood estimator for  cooperative sequential adsorption model  (CSA) 
were studied. CSA is a probabilistic model motivated by 
adsorption processes in physics and chemistry (\cite{Evans}).
The main peculiarity of  adsorption processes  is that adsorbed particles  
change adsorption properties  of the material. For instance, the subsequent 
particles might be more likely to be adsorbed around locations of previously 
adsorbed particles. In other words, the adsorption process might
accelerate as the surface gets saturated.
In the opposite scenario adsorbed particles inhibit adsorption of subsequent particles, so that
 the adsorption process slows down. 

Mathematically CSA is formulated as a random  sequential allocation of points 
in a bounded region of space (the observation window). 
The result of CSA dynamics is a sequential point pattern
which seems to be of great interest in many applications. 
It should be noted that CSA  can produce  a large variety of aggregated point 
patterns (see, e.g., the images throughout the paper).  

It was first  noticed by physicists 
(e.g., see \cite{Evans}, p.1285) that this type of  model 
can be used for   modelling the spatial-temporal processes
 similar to the irreversible 
spread of disease or epidemics. This idea is developed further
 in  \cite{MPVS}
where use of CSA for  modelling time series of spatial locations is discussed.

Biological growth was mentioned in \cite{Evans} as another  potential application 
of the adsorption models. 
These ideas have been recently supported by both experimental and simulation studies 
of keratin filament (KF)  network formation in biology.
KF networks are part of the cell cytoskeleton 
and they  determine the shape and biophysical properties of the cells. 
 Loosely speaking, the KF is an aggregated spatial structure  
formed by a union of curved finite segments (fibres). 
Experimental results  (\cite{Windoffer}) and 
simulation studies  (\cite{Beil}) suggest that  the KF can be thought 
 as a result of a  sequential spatial growth process with self-organising properties.
It is also argued  in \cite{Marshall} (see also references therein), that
self-organizing processes combined with simple physical constraints seem to
 have key roles in controlling organelle size, number, shape and position, 
and these factors then combine to produce the overall cell architecture.  
CSA  seems to be useful for modelling  spatial random growth with self-organising  properties. 

The variant of  CSA under consideration here is easy to parameterise. 
Statistical inference for the model parameters developed in \cite{MPVS}
was based on maximum  likelihood estimation (MLE).
It was shown in \cite{MPVS} that maximum likelihood estimator exists uniquely. Moreover,
it was proved that the maximum likelihood estimator  is  consistent in the 
thermodynamic limit.
The thermodynamic limit means that the observation window expands
to the whole space and the number of allocated points grows linearly in the 
volume of the window. 
The main result of the present paper 
is asymptotic normality of maximum likelihood estimator in the same limit.

The study of statistical properties
 of MLE  in both \cite{MPVS} and  this paper 
 is essentially based on the fact that the model likelihood
depends on the point configuration via  
statistics with  a certain special structure,
allowing us  to apply   the limit theory  for random sequential packing 
and deposition  (see, e.g., \cite{PenYuk}).

\section{CSA as a generalisation of random sequential adsorption}

The adsorption model most commonly
 studied in the physics literature is random sequential adsorption (RSA).
Mathematically RSA is formulated  as the following  packing model. 
Consider  a bounded region $D$  of Euclidean space
$\R^d$ (modelling the adsorbing material) and a sequence of  
 independent  points
$Y_1, Y_2, \ldots,$ (modelling 
the particles) sequentially arriving in $D$ at random. 
 An arriving point  is accepted with probability  $1$, 
if the ball of a certain fixed radius $R$ (interaction radius)
centered  at the point  does not cover  any of previously 
accepted (adsorbed) points;
otherwise the point is rejected.

\begin{figure}[htbp]
\centering
\begin{tabular}{cc}
    \includegraphics[width=3.2in, height=2.9in, angle=270]{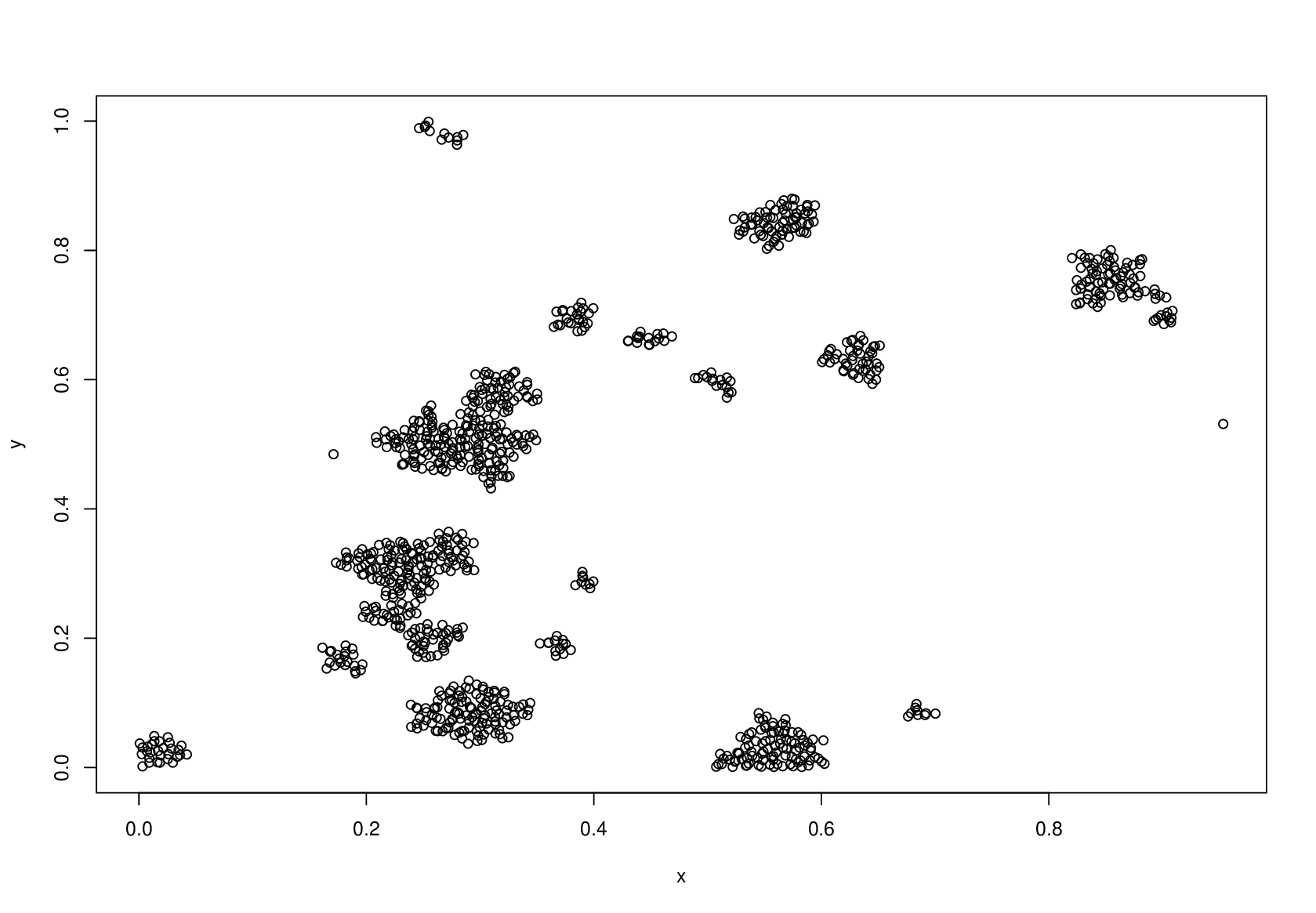}&
  \includegraphics[width=3.2in, height=2.4in, angle=270]{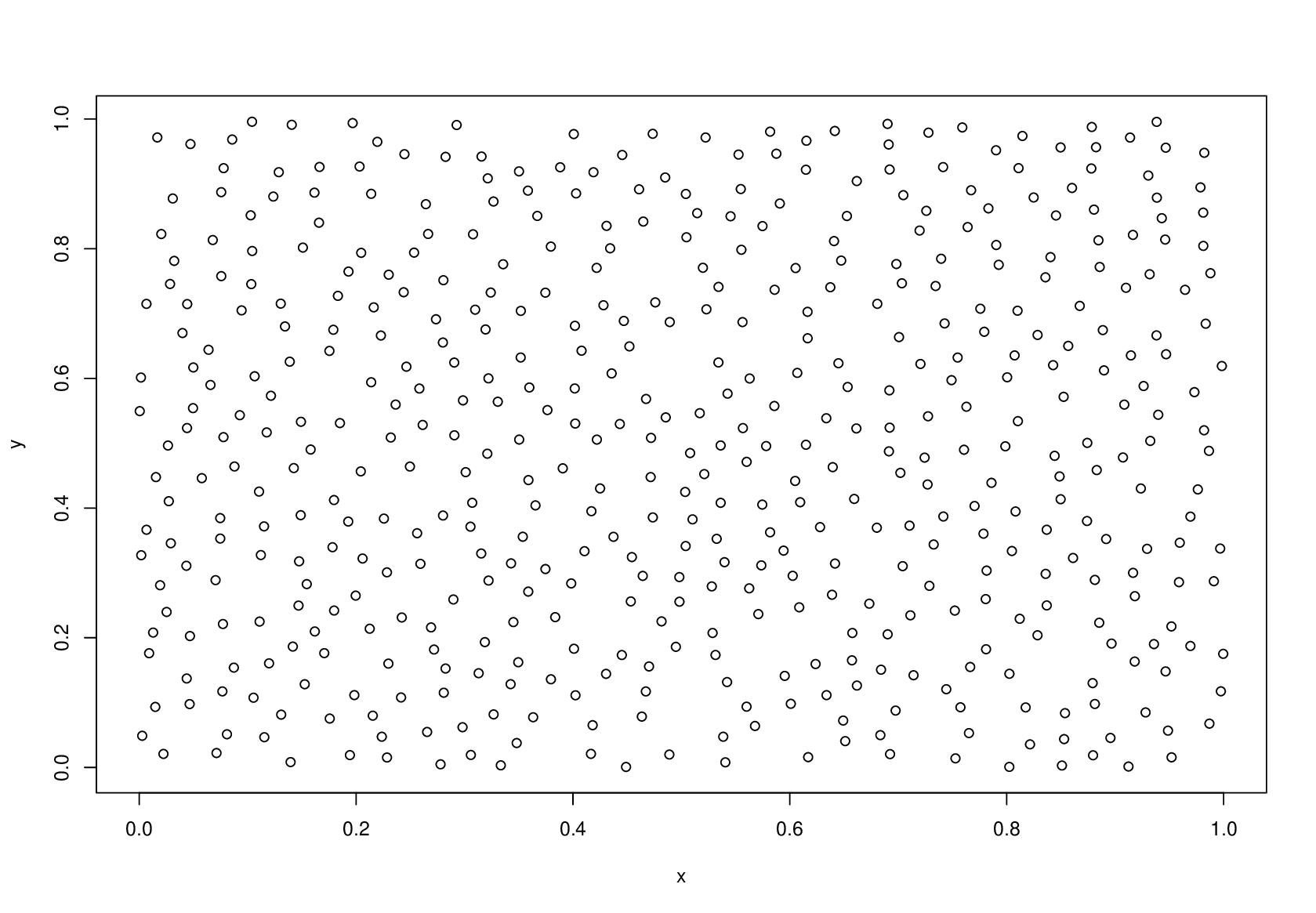}
\end{tabular}
\label{fig1}
\end{figure}

RSA with interaction radius $R$ is nothing else but  the
$d$-dimensional version of the classical
  car parking model \cite{Renyi}, where a ``car'' is
 modelled by a ball of radius $R/2$.
Clearly the distance between any two points in a RSA point pattern 
cannot be less  than the interaction radius $R$. 
Therefore RSA  generates   only  regular point patterns which
 are similar to the right  one in Figure \ref{fig1}, 
and never generates  point patterns  similar to the left one  in 
Figure \ref{fig1}.
However, RSA can  be easily  generalised  in order to generate
  aggregated point patterns.
To do so, we  allow {\it neighbours}. That is, we
let an arriving  point be   accepted with a certain conditional probability,
even if a ball of radius $R$ 
centred at the point covers some of the previously accepted points. 
In  general, the acceptance probability can depend     
on the spatial configuration formed by previously accepted points.
We study the model in which    the acceptance probability depends on the number of neighbours. 

More precisely,  fix a sequence of  non-negative  numbers $\beta_0, \beta_1,...$, such that 
$\beta_0>0$.
Given a  sequence of  accepted points $X(k)=(X_1,\ldots, X_k)$ ($X(0)=\emptyset$),
let the next uniform arrival $Y$   
be accepted with conditional probability proportional to $\beta_i$, 
if the number of neighbours of $Y$ among $X_1,\ldots,X_k$ is equal to $i\geq 0$.
If $\beta_0>0$ and $\beta_k=0,\, k\geq 1$, then this model is RSA.

This CSA model can be regarded as  a continuous version  
 of the  lattice model (i.e. where $D$ is a subset of lattice $\Z^d$) 
known as monomer filling with nearest-neighbour cooperative effects. 
CSA in this particular form was formulated for the first time
 in \cite{Sch1}, where 
its asymptotic study  was undertaken under certain assumptions.
 In what follows we denote by CSA the adsorption model of this type.

CSA can be used   for modelling both clustered and regular point patterns.
A large variety of aggregated point patterns can be generated 
by modulating the model parameters. 
For instance,  the left image  in Figure \ref{fig1}, containing $1000$ points,  is generated  by CSA with parameters 
$R=0.01, \beta_0=1, \beta_1=1000, \beta_2=10000, \beta_k=0,\, k\geq 3$. The right image
(containing $500$ points) is a typical regular image produced by RSA 
(here the interaction radius is $R=0.03$).

\section{Notation and assumptions}

Let $D$ be a convex compact subset of $\R^d$, $R$ be a positive constant, and  
$\{\beta_k,\, k\geq 0\}$ be a sequence of non-negative numbers.
For any point $x\in \R^d$ and any finite sequence
$\by=(y_1,\ldots,y_n),\, n\geq 1,$ of points in 
$\R^d$,
 we denote by $\nn(x,\by)$ the number of points $y_i $ in the sequence
$\by$,
such that the distance between $x$ and $y_i$ is not greater than
$R$.
By definition $\nn(x,\emptyset)=0.$

Let  $X(\ell)=(X_1,\ldots,X_{\ell}),\, X_i\in \R^d,\, i=1,\ldots,\ell$
be a vector of first $\ell$ random points  
sequentially generated  by CSA.
CSA dynamics goes as follows.
Given a sequence of accepted points  $X(k)=(X_1,\ldots,X_{k})$ (which can
 be empty, i.e. $k=0$) 
a new point  $Y$, uniformly distributed in $D$,
 is accepted with probability proportional to $\beta_{\nn(Y, X(k))}$
and rejected otherwise.
If $Y$ is accepted, then we set $X_{k+1}=Y$ and $X(k+1)=(X_1,\ldots,X_k, X_{k+1})$.
The conditional 
probability density function of the next accepted point $X_{k+1}$
is 
\begin{equation}
\label{dens}
\psi_{k+1}(x)=\frac{\beta_{\nn(x, X(k))}}{\int_{D}\beta_{\nn(y, X(k))}dy},\,\, x\in D.
\end{equation}
It is easy to see that 
the sequence of accepted points is an embedded Markov chain for a 
continuous time spatial birth process $\bx(t)\in D,\,t\geq 0,$  specified 
by the following birth rates. If the process state 
at time $t\geq 0$ is $\bx$, then  the birth rate at point  $x\in D$ is  
$\beta_{\nn(x,\bx)},$ 
the total birth rate is  
$$\alpha(\bx)=\int\limits_{D} \beta_{\nn(x,\bx)}dx,$$ 
and the waiting time until the next process jump is an 
exponential random variable with mean  $\alpha^{-1}(\bx)$.

As in \cite{MPVS}, we assume throughout that 
\begin{itemize}
\item there is a {\it finite} number of positive $\beta's$, that is 
$\beta_0>0,\,\ldots, \beta_N>0$
and  $\beta_k=0,$ for $k\geq N+1$, for some $N\geq 1$, where the number $N$ can be  {\it unknown},
\item $\beta_0=1$,
\item the interaction radius $R$ is a {\it fixed} and  {\it known} constant. 
\end{itemize}
It is easy to see that the  joint probability density
$\prod_{k=1}^\ell \psi_k(x_k)$ of the first $\ell$ accepted 
points  can be written as follows:
\begin{equation}
\label{like}
p_{ \ell, \beta,D}(x_1,\ldots,x_\ell)=\frac{\prod_{k=0}^{N}\beta_{k}^{t_{k}(x(\ell))}}
{\prod_{k=1}^{\ell}\int_{D}\beta_{\nn(x, x(k-1))}dx}
{\bf 1}_{\{\widehat{N}(x(\ell))\leq N\}},
\end{equation}
where
\begin{equation}
\label{Nestim} 
\widehat{N}(x(\ell))=\max\limits_{x_i\in x(\ell)}\nn(x_i, x(i-1)),
\end{equation}
and  
\bea
t_k(x(\ell))=
\sum\limits_{i =1}^\ell
 {\bf 1}_{\{\nn(x_i,x(i-1))=k\}},\,\, 
k=0,\ldots,N.
\label{tkdef}
\eea
where we denoted for short 
$x(k)=(x_1,\ldots,x_{k}),\, k\geq 1,$
and $x(0)=\emptyset$ for $k=0$.

{\bf Remark.}
{\it 
It should be noticed that we do not completely recover the parameters of the spatial birth process. 
In the present setting 
we do statistical inference only for  the embedded Markov chain, which 
distribution is completely specified by the ratios $\beta_i/\beta_0,\,
1,\ldots,N$ and the interaction radius. As a result, 
one can forecast the probability  distribution of the next accepted point, but 
not the waiting time  until  the  next acceptance event.}

\bigskip 

As in \cite{MPVS}, let $D_1$ be  the unit cube centred at the origin and 
consider a sequence of rescaled domains 
$$D_m=m^{1/d}D_1,\, m\in \Z_{+}.$$
Fix  $\{\ell_m,\, m\geq 1\}$ an arbitrary monotonically increasing sequence of 
positive numbers, where $\ell_m$ stands for the number of observed points 
in the domain $D_m$.
\begin{assum}
\label{A1}
The number of observed points
is  asymptotically linear in $m$, that is 
$$
\lim_{m \to \infty}  \left( \frac{\ell_m}{m}
\right)
=  \mu \in (0,\theta_\infty),$$
where $\theta_\infty$ is the jamming density (\cite{MPVS}).
\end{assum}
Define
$$\SSS_m: =\{x(\ell_m)=(x_1,\ldots,x_{\ell_m}),\,x_i\in D_m: 
\widehat{N}(x(\ell_m))\leq N\}.$$
Given parameters  $\beta=(\beta_1,\ldots,\beta_N)$
consider  a probability measure $\P_{m, \beta}$  on $\SSS_m$
specified by the probability density (\ref{like}) with $\ell =\ell_m$
and $D = D_m$. 
Expectation with respect to this measure
is denoted by $\E_{m, \beta}$.
We assume that $\beta\in {\cal B}$, where ${\cal B}$ is an open subset of $\R^N$, such that ${\cal B}\subset \R_{+}^N$.
The true parameter is denoted by  $\beta^{\left(0\right)}=
\left(\beta_1^{\left(0\right)},\ldots,\beta_N^{\left(0\right)}\right)$. 
Also, we denote for short $\P_m^{\left(0\right)}=\P_{m, \beta^{\left(0\right)}}$
 and $\E_m^{\left(0\right)}=\E_{m, \beta^{\left(0\right)}}$.

\section{The results} 
\label{results}

Given $m$ assume  $\ell_m \geq 2$ and consider the log likelihood function
\begin{equation}
\label{log0}
L_m(X^m(\ell_m), \beta) = \log(p_{\ell_m, \beta,D_m}(X^m_1,\ldots,X^m_{\ell_m})),
\end{equation}
where $X^{m}(\ell_m)$ is the vector of observed points  in $D_m$. 
Given observation $X^{m}(\ell_m)$ we define the maximum likelihood estimators
$$
\widehat{\beta}(X^{m}(\ell_m))=(\widehat{\beta}_{1,m},\ldots,\widehat{\beta}_{N,m})$$ 
of parameters $\beta^{\left(0\right)}=(\beta_1^{\left(0\right)},\ldots,
\beta_{N}^{\left(0\right)})$  as 
maximizers of function $L_m(X^{m}(\ell_m), \beta)$ and  which 
can be found   as a solution of the following  system of MLE equations 
\begin{equation}
\label{mle0}
\frac{\partial L_m(X^{m}(\ell_m), \beta)}{\partial \beta_j}=0,\quad j=1,\ldots,N.
\end{equation}

The following two statements were proved in \cite{MPVS} 
(see Theorem $2.2$ and Lemma $5.2$, part $2)$,
respectively in \cite{MPVS}).
\begin{lm}
\label{L1}
Under Assumption \ref{A1}  with $\P_m^{\left(0\right)}-$probability
tending to $1$ as $m\ri \infty$  there exists a unique  
positive solution $(\widehat{\beta}_{1,m},\ldots,\widehat{\beta}_{N,m})$ of the likelihood 
equations and 
$$ (\widehat{\beta}_{1,m},\ldots,\widehat{\beta}_{N,m}) \to
(\beta_1^{\left(0\right)},\ldots,\beta_N^{\left(0\right)})$$
in $\P_m^{\left(0\right)}-$ probability as $m\ri \infty$.
\end{lm}
\begin{lm}
\label{L2}
Consider the   matrix 
$$J_{m}(X^{m}(\ell_m), \beta) :=-\left(\frac{\partial^2  L_m(X^{m}(\ell_m), \beta)}
{\partial \beta_i \partial \beta_j} \right)_{i,j=1}^{N}.$$
There is a family of $N \times N$ real matrices $J(\beta,\mu)$, defined
for $\beta \in {\cal B}$ and $\mu \in(0,\theta_\infty)$,
 such that under Assumption \ref{A1} 
$$-\frac{J_{m}(X^{m}(\ell_m), \beta)}{m}\ri J(\beta, \mu)$$
in $\P_m^{\left(0\right)}-$ probability as $m\ri \infty$ for
 any $\beta\in {\cal B}$. Moreover,
the limit matrix evaluated at 
$\beta=\beta^{\left(0\right)}$, i.e.
\begin{equation}
\label{info}
J^{\left(0\right)}(\mu)=J\left(\beta^{\left(0\right)}, \mu\right)
\end{equation}
is positive definite. Finally, if $\beta(m)$ is a {\em random}
${\cal B}$-valued sequence converging in probability to
$\beta^{(0)}$ as $m \to \infty$, then
$$-\frac{J_{m}(X^{m}(\ell_m), \beta(m))}{m}\ri J(\beta^{(0)}, \mu)$$
in $\P_m^{\left(0\right)}-$ probability as $m\ri \infty$.
\end{lm} 
The last part of Lemma \ref{L2} is not included in
Lemma 5.2(2) of \cite{MPVS}, but can be proved in the same
manner as that result.

In Section \ref{matr} we give extended study  of  the structure of 
 the limit information matrix.

\begin{thr}
\label{L3}
Under Assumption \ref{A1} the model score function 
\begin{equation}
\label{score}
\nabla L_m(X^{m}(\ell_m), \beta^{(0)})=
\left.
\left(\frac{\partial L_m(X^{m}(\ell_m), \beta)}{\partial 
\beta_1}, \ldots, \frac{\partial L_m(X^{m}(\ell_m), \beta)}{\partial 
\beta_N}\right) \right|_{\beta = \beta^{(0)}}
\end{equation}
converges in distribution as $m\ri \infty$
 to a Gaussian vector with mean zero  and 
 covariance matrix $J^{\left(0\right)}(\mu)$.
\end{thr}
Theorem  \ref{L3} is proved in Section \ref{CLT}.
The following theorem states that the
 MLE is asymptotically normal. This is the main  result of the paper.
\begin{thr}
\label{main}
Under Assumption \ref{A1}
$$\sqrt{m}\left(\widehat{\beta}(X^{m}(\ell_m))-\beta^{\left(0\right)}\right)\ri 
{\cal  N}\left(0, \left(J^{\left(0\right)}(\mu)\right)^{-1}\right)$$
in distribution as $m\ri \infty$, where  
${\cal N}\left(0, \left(J^{\left(0\right)}\right)^{-1}(\mu)\right)$ 
is the Gaussian vector with zero mean and  with the covariance matrix 
$\left(J^{\left(0\right)}(\mu)\right)^{-1}$.
\end{thr}
Theorem \ref{main} provides asymptotic justification for  creating
 confidence intervals based on the normal distribution,
as we do in the example in Section \ref{example}. 

\section{The model likelihood}

In this section we introduce  more notation and recall some
other facts from \cite{MPVS} which will be used in Section \ref{proofs}.

Let $X^m(\ell_m)=(X^m_1,\ldots,X^m_{\ell_m})$ be the sequence of 
observed  points $X_i^m$ in $ D_m$. 
 Denote
\begin{equation}
t^m_{j,k}=
t_{j}(X^m(k)) \quad 0\leq k\leq \ell_m-1,\, j=1,\ldots,N,
\label{tjmdef}
\end{equation}
where $t_j,\, j=1,\ldots,N,$ are statistics defined by  equation (\ref{tkdef}), and denote 
\begin{equation}
\label{Gajmdef}
\Gamma_{j,k}^m=\Gamma_{j, k}^m(X^m(k))
=\int\limits_{D_m}{\bf 1}_{\{u:\nn(u,X^m(k))=j\}} du, \quad 0\leq k\leq \ell_m-1,\, j\geq 0,
\end{equation}
note that $\Gamma_{j,k}^m(X^m(k))=0$ for $k<j$ and that $\Gamma_{0,0}^m(X^m(k))$  
is equal to $m$.

In terms of $t-$ and $\Gamma-$statistics, using (\ref{like}), (\ref{tkdef})
and (\ref{Gajmdef})
  the model likelihood can be rewritten as follows 
\begin{align}
\nonumber
L_m(X^m(\ell_m), \beta) &= \log(p_{\ell_m, \beta,D_m}(X^m_1,\ldots,X^m_{\ell_m}))\\
&=\sum_{k=1}^{N}t_{k}(X^{m}(\ell_m))\log(\beta_k)
-\sum\limits_{k=1}^{\ell_m}\log\left(\int_{D_m} \beta_{\nn(u,X^m(k-1))}
 du\right) \nonumber\\
&=\sum_{k=1}^{N}t_{k, \ell_m}^m\log(\beta_k)
-\sum\limits_{k=1}^{\ell_m}\log\left(\Gamma_{0,k-1}^m + 
\sum\limits_{j=1}^{N}\beta_{j}\Gamma_{j, k-1}^m\right).\label{0925a}
\end{align}

Thus  the log likelihood function  depends on the 
observed point configuration only through  $t-$statistics 
$t_{k}$ and $\Gamma-$statistics  $\Gamma_{j, k}^m$. 

Theorem 2.2 in \cite{MPVS} says that 
if 
$$\lim_{m \to \infty} (\ell_m/m) = \mu \in (0,\theta_\infty(\beta)),$$
then  as $m \to \infty$ we have for any $\beta \in {\cal B}$ that 
 \bea
\frac{t^m_{j,\ell_m}}{m}\toP \rho_j\left(\mu, \beta\right),\, j=1,\ldots,N,
\label{1214a}
\eea
and 
\bea 
\frac{\Gamma^m_{j,\ell_m}}{m}\toP \gamma_{j}\left(\mu, \beta\right),\, j=0,\ldots,N,
\label{1214b}
\eea 
where the functions 
$(\rho_j\left(\mu, \beta\right), 
\mu \in (0,\theta_\infty(\beta)),1 \leq j \leq N
$
and $(\gamma_{j}\left(\mu, \beta\right), 
\mu \in (0,\theta_\infty(\beta)), 0 \leq j \leq N
$
are  strictly positive and continuous in $\mu$, and are
 related by the following integral equation
\begin{equation}
\label{t_g}
\rho_j\left(\mu, \beta\right)=\int\limits_{0}^{\mu}
\frac{\beta_j\gamma_{j}\left(\lambda, \beta\right)}{\gamma_{0}\left(\lambda, \beta\right)+
\sum_{i=1}^N\beta_i \gamma_{i}\left(\lambda, \beta\right)}d\lambda,\, j=1,\ldots,N,
\end{equation}
for any  $0 < \mu < \theta_{\infty}(\beta)$. 

\section{Proofs}
\label{proofs}

\subsection{Proof of Theorem \ref{main}}

Given Theorem \ref{L3},
the proof of Theorem \ref{main}, although new to this particular model,
 runs along standard lines
(see e.g.  \cite{Leh}, or Theorem 1 of \cite{BRR}),
and we give just a sketch.

Choose $\delta >0$ such that the ball of radius $\delta$
centered at $\beta^{(0)}$ is contained in ${\cal B}$.
By consistency of the maximum likelihood estimator
 $\widehat{\beta}(X^m(\ell_m))$ (Lemma \ref{L1}), 
we have that 
$$|\widehat{\beta}(X^m(\ell_m))-\beta^{\left(0\right)}|<\delta,$$
 with probability $\P_m^{\left(0\right)}$ close to $1$ if $m$ is 
large enough.
With $\partial_i$ denoting differentiation with respect to the
$i$th component of $\beta$,
we make a Taylor expansion of 
$\partial_i  (L_m(X^m(\ell_m), \beta)$ about $\beta^{\left(0\right)}$:
\begin{align*}
0= \partial_i
  L_m(X^m(\ell_m), \widehat{\beta}(X^m(\ell_m)))&=
\partial_i   L_m(X^m(\ell_m),\beta^{\left(0\right)})\\
&+ 
\sum_{j=1}^N \partial^2_{ij} (X^m(\ell_m),\bar{\beta})
(\widehat{\beta}
(X^m(\ell_m))-\beta^{\left(0\right)})_j,
\end{align*}
where $\bar{\beta}$ lies on the line segment from $\beta^{\left(0\right)}$
to $\widehat{\beta}(X^m(\ell_m))$.
Rewriting this expression, we obtain
$$
\sum_{j=1}^N \left( \frac{ 
 -\partial_{ij}^2(X^m(\ell_m),\bar{\beta})}{m} \right)
\left( 
\sqrt{m}
(\widehat{\beta}(X^m(\ell_m))-\beta^{\left(0\right)})_j \right)=
\frac{\partial_i
  L_m(X^m(\ell_m),\beta^{\left(0\right)})}{\sqrt{m}}.$$
In the left hand expression $\bar{\beta}$ depends on $i$
but converges in probability to $\beta^{(0)}$
as $n \to \infty$ by  Lemma \ref{L1}. By Lemma \ref{L2},  
for each $(i,j)$ the first factor inside the sum
converges in probability to $J_{ij}^{(0)}(\mu)$.
Observing that
Theorem  \ref{L3} 
applies to the right hand side,
we can complete 
the proof by applying Lemma 6.4.1 of
\cite{Leh}.

\subsection{Proof of Theorem  \ref{L3}}
\label{CLT}

Let ${\cal F}_{j}^{\left(m\right)}=\sigma\{X^m_1,\ldots,X^m_{j}\}$ be the 
$\sigma-$algebra generated by   the first $j$ points  observed in $D_m$.
Asymptotic normality of the score function is  essentially based  on the
 following fact. Namely, for any $k=1,\ldots,N$, the triangle array 
\begin{equation}
\label{trian}
\left\{\E_{m, \beta}\left(\left.\frac{\partial L_m(X^m(\ell_m), \beta)}{\partial \beta_k}\right|
 {\cal F}_{j}^{\left(m\right)}\right), 
{\cal F}_{j}^{\left(m\right)}\right\}_{j=1}^{\ell_m},\, m\geq 2,
\end{equation}
 is a zero-mean square integrable martingale array.
Indeed, by the representation (\ref{0925a}),  
\begin{equation}
\label{score1}
\frac{\partial L_m(X^m(\ell_m), \beta)}{\partial 
\beta_k}=\frac{t_k(X^m(\ell_m))}{\beta_k} - \sum\limits_{j=1}^{\ell_m}
 \frac{\Gamma^m_{k,j-1}}{\Gamma^m_{0,j-1} + 
\sum_{i=1}^{N}\beta_{i}\Gamma^m_{i, j-1}},
\end{equation}
for $j=1,\ldots, N.$
Introducing  the following quantities
\begin{equation}
\label{xi}
\xi^{m}_{k,i}={\bf 1}_{\{\nn(X^m_i,X^m(i-1))=k\}},\,\, 
k=0,\ldots,N,\,\, i=1,\ldots,\ell_m,
\end{equation}
allows to rewrite equation (\ref{tkdef})  as follows:
$$
t_k(X^m(\ell_m))=
\sum\limits_{i =1}^{\ell_m}
\xi^{m}_{k,i},\,\, 
k=0,\ldots,N.
$$
Denote for short 
$$\bar{\xi}^{\, m}_{k,i}=\E_{m, \beta}\left(\xi^{m}_{k,i}|
{\cal F}_{i-1}^{\left(m\right)}\right).$$
It is easy to see that 
\begin{equation}
\bar{\xi}^{\, m}_{k,i}=\frac{\beta_k\Gamma^m_{k,i-1}}{\Gamma^m_{0,i-1} + 
\sum_{j=1}^{N}\beta_{j}\Gamma^m_{j, i-1}}
,\,\, 
k=1,\ldots,N,\,\, i=1,\ldots,\ell_m.
\label{100504a}
\end{equation}
By using notation
\begin{equation}
\label{zeta}
\zeta^m_{k, i}=\frac{1}{\beta_k}
\left(\xi^m_{k, i}-
\bar{\xi}^{\, m}_{k, i}\right)
\end{equation}
equation  (\ref{score1}) can now be rewritten as follows:
\begin{equation}
\label{score2}
\frac{\partial L_m(X^m(\ell_m), \beta)}{\partial 
\beta_k}=\frac{1}{\beta_k}\sum\limits_{i=1}^{l_m}\left(\xi^m_{k,i}-
\bar{\xi}_{k,i}^{\, m}\right)=
\sum\limits_{i=1}^{l_m}\zeta^m_{k,i},
\quad k=1,\ldots, N.
\end{equation}
Therefore the triangle array (\ref{trian}) 
 is a zero-mean square integrable martingale array
with differences given by equation (\ref{zeta}). 
This implies that for any real vector ${\bf a}=(a_1,\ldots, a_N)^T$
$$\left\{ \E_{m, \beta}\left(\left. \frac{1}{\sqrt{m}}\sum\limits_{k=1}^N a_k
\frac{\partial L_m(X^m(\ell_m), \beta)}{\partial \beta_k}\right|{\cal F}_{j}^{\left(m\right)}\right),\,
{\cal F}_{j}^{\left(m\right)}\right\}_{j=1}^{\ell_m},\, m\geq 2,$$
 is a zero-mean square integrable martingale array.

By the Cram\'er-Wold device (see for example \cite{Bill}),
Theorem  \ref{L3}   follows from
the following fact.
 
\begin{lm}
\label{univar}
Under Assumption \ref{A1} for any real vector ${\bf a}=(a_1,\ldots, a_N)^T$,
$$\frac{1}{\sqrt{m}}\sum\limits_{k=1}^N a_k
\frac{\partial L_m(X^m(\ell_m), \beta^{\left(0\right)})}{\partial \beta_k} \ri 
{\cal N}\left(0, \sigma^2_{\bf a}\right),$$
in distribution as $m\ri \infty$, where 
$$\sigma^2_{\bf a}={\bf a}^TJ^{\left(0\right)}(\mu){\bf a},$$
$J^{\left(0\right)}(\mu)$ is the matrix defined by equation (\ref{info})
and 
${\cal N}\left(0, \sigma^2_{\bf a}\right)$
is the Gaussian vector with zero mean and variance $\sigma^2_a$.
\end{lm}

In proving Lemma \ref{univar} we shall repeatedly use 
the following fact which is simple enough for us to omit its proof.
\begin{prop}
\label{dct_prob}
Let $\xi_n,\,n\geq 1,$ and $\eta_n,\, n\geq 1,$ be two  sequences of random variables, $C>0, a$ and $b$ be some constants.
Suppose that  $|\xi_n|<C, |\eta_n|<C$, $\xi_n \ri a$ in probability as $n\ri \infty$
and $\E(\eta_n)\ri b$ as $n\ri \infty$.
Then 
$\E(\xi_n\eta_n)\ri ab$ as $n\ri \infty$.
\end{prop}

{\em Proof of Lemma \ref{univar}.}
By (\ref{score2}),
 for any $\beta \in {\cal B}$
we have
$$\frac{1}{\sqrt{m}}\sum\limits_{k=1}^N a_k
\frac{\partial L_m(X^m(\ell_m), \beta)}{\partial \beta_k}=\frac{1}{\sqrt{m}}
\sum\limits_{i=2}^{\ell_m}\eta^m_i$$
 where
\begin{equation}
\label{zeta1}
\eta^m_i=\sum\limits_{k=1}^N a_k\zeta^m_{k, i}
\end{equation}
and $\zeta^m_{k, i}$ are the quantities defined by equation (\ref{zeta}).
It is easy to see that
\begin{equation}
\label{i}
\frac{1}{\sqrt{m}}\max\limits_{i}|\eta^m_{i}|\leq \frac{2N}{\sqrt{m}}\max\limits_{k=1,\ldots,N}
\left(\frac{a_k}{\beta_k}\right)
\ri 0, \quad \mbox{as} \quad m\ri \infty,
\end{equation}
and
\begin{equation}
\label{ii}
\frac{1}{m}\E_{m, \beta} \left(\max\limits_{i}\left(\eta_i^{m}\right)^2\right)\leq \frac{4N^2}{m}
\max\limits_{k=1,\ldots,N}
\left(\frac{a_k}{\beta_k}\right)^2
\ri 0,\quad \mbox{as} \quad m\ri \infty.
\end{equation}
By Propositions \ref{prop1} and \ref{prop2} below, 
we also have
 under Assumption \ref{A1} that
\begin{equation}
\label{2nd_cond}
\frac{1}{m}\sum\limits_{i=2}^{\ell_m}
(\eta_{i}^m)^2 \ri {\bf a}^TJ^{\left(0\right)}(\mu){\bf a}
\end{equation}
in $\P^{(0)}_{m}-$probability as $m\ri \infty$.

Using (\ref{ii}), (\ref{ii}) and (\ref{2nd_cond}),
we can then 
 apply the central limit theorem for martingale difference arrays
(Theorem (2.3) of \cite{McLeish}) to complete the proof of
Lemma \ref{univar}.

\begin{prop}
\label{prop1}
Under Assumption \ref{A1}
\begin{equation}
\label{1st_cond}
\lim\limits_{m\ri \infty}\frac{1}{m}\sum\limits_{i=2}^{\ell_m}
\E_m^{\left(0\right)}\left((\eta^m_{i})^2\right)
={\bf a}^TJ^{\left(0\right)}(\mu){\bf a}.
\end{equation}
\end{prop}
{\it Proof of Proposition \ref{prop1}.}
It was shown in Section 6.2 of  \cite{MPVS} that
 the limit of the scaled Hessian in Lemma \ref{L2} 
evaluated at the true parameter has the following integral representation
\begin{equation}
\label{integ}
J^{\left(0\right)}(\mu)=J\left(\beta^{\left(0\right)}, \mu\right)
=\int\limits_{0}^{\mu}Q^{\left(0\right)}(\lambda)d\lambda,
\end{equation}
where
\begin{equation}
\label{mat}
Q^{\left(0\right)}(\lambda)=
\left(\frac{ \gamma_i^{(0)}(\lambda) }{\beta_i^{(0)} Z(\beta^{\left(0\right)},\lambda)}
 \delta_{ij}
- \frac{ \gamma_i^{(0)}(\lambda) \gamma_j^{(0)}(\lambda) }{Z^2(\beta^{\left(0\right)},\lambda)}\right)_{i,j=1}^N,  
\end{equation}
where $\delta_{ij}$ is the Kroneker symbol,
$\gamma_{j}^{\left(0\right)}(\lambda)=
\gamma_{j}(\lambda, \beta^{\left(0\right)}),\, j=0,\ldots,N$ ($\gamma-$functions 
are defined by (\ref{1214b})) 
and 
\bea
Z(\beta, \lambda)=\gamma_0^{\left(0\right)}(\lambda)+
\sum_{i=1}^N\beta_i \gamma_{i}^{\left(0\right)}(\lambda).
\label{071218a}
\eea
Let us show that if $i=i_m$ is such that $i/m\ri \lambda \in (0,\mu)$, as $m\ri \infty$,
then 
$$\E_m^{\left(0\right)}\left((\eta_i^m)^2\right)\ri 
{\bf a}^TQ^{\left(0\right)}(\lambda){\bf a}$$
as $m\ri \infty$.
Indeed, 
\begin{align*}
\E_m^{\left(0\right)}\left((\eta_i^m)^2\right)&=\sum\limits_{k,j=1}^N a_ka_j\E_m^{\left(0\right)}(\zeta^m_{k,i}\zeta^m_{j,i})\\
&=\sum\limits_{k,j=1}^N 
\frac{a_ka_j}{\beta^{\left(0\right)}_k\beta^{\left(0\right)}_j}
\E_m^{\left(0\right)}((\xi^m_{k,i}-\bar{\xi}^{\, m}_{k,i})(\xi^m_{j,i}-\bar{\xi}^{\, m}_{j,i}))\\
&=\sum\limits_{k,j=1}^N \frac{a_ka_j}{\beta^{\left(0\right)}_k\beta^{\left(0\right)}_j}
\E_m^{\left(0\right)}
\left(\xi_{k,i}^m\xi_{j,i}^m-\xi_{k,i}^m\bar{\xi}^{\,m}_{j,i}-
\xi_{j,i}^m\bar{\xi}^{\, m}_{k,i}
+\bar{\xi}_{k,i}^{\, m}\bar{\xi}^{\, m}_{j,i}\right)
\end{align*}
Notice that 
$$\E_m^{\left(0\right)}(\xi_{k,i}^m\xi_{j,i}^m)=
\E_m^{\left(0\right)}(\xi_{k,i}^m)\delta_{kj}=
\E_m^{\left(0\right)}(\bar{\xi}_{k,i}^{\, m})\delta_{kj},$$
where $\delta_{ij}$ is the Kroneker symbol and 
$$\E_m^{\left(0\right)}(\xi_{k,i}^m\bar{\xi}_{j,i}^{\, m})=
\E_m^{\left(0\right)}(\xi_{j,i}^m\bar{\xi}_{k,i}^{\, m})=
\E_m^{\left(0\right)}(\bar{\xi}_{k,i}^{\, m}\bar{\xi}_{j,i}^{\, m}).$$
Therefore
$$
\E_m^{\left(0\right)}\left((\eta_i^m)^2\right)=\sum\limits_{k,j=1}^N \frac{a_ka_j}{\beta^{\left(0\right)}_k\beta^{\left(0\right)}_j}
\left[\E_m^{\left(0\right)}(\bar{\xi}_{k,i}^{\, m})\delta_{kj}-
\E_m^{\left(0\right)}(\bar{\xi}_{k,i}^{\, m}\bar{\xi}_{j,i}^{\, m})\right].
$$

 By (\ref{100504a}) and (\ref{1214b}),
\begin{equation}
\label{xi_conv}
\bar{\xi}^{\, m}_{r,i}=
\frac{\beta_r^{\left(0\right)}\Gamma^m_{r,i-1}}{\Gamma^m_{0,i-1} + 
\sum_{j=1}^{N}\beta_{j}^{\left(0\right)}\Gamma^m_{j, i-1}}
 \ri
 \frac{\beta_r^{\left(0\right)}\gamma_{r}^{\left(0\right)}(\lambda)}
{\gamma_0^{\left(0\right)}(\lambda)+\sum_{j=1}^N\beta_j^{\left(0\right)} \gamma_{j}^{\left(0\right)}(\lambda)}
=
\frac{\beta_r^{\left(0\right)}\gamma_{r}^{\left(0\right)}(\lambda)}{Z\left(\beta^{\left(0\right)}, \lambda\right)}
\end{equation}
in $\P^{(0)}_{m}$ probability as $i/m\ri \lambda$, for any $r=0,\ldots,N$.
This fact along with Proposition \ref{dct_prob}
yield that 
\begin{multline*}\frac{1}{\beta^{\left(0\right)}_k\beta^{\left(0\right)}_j}
\left[\E_m^{\left(0\right)}(\bar{\xi}_{k,i}^{\, m})\delta_{kj}-
\E_m^{\left(0\right)}(\bar{\xi}_{k,i}^{\, m}\bar{\xi}_{j,i}^{\, m})\right]\\
\ri 
\frac{ \gamma_k^{(0)}(\lambda) }{\beta_k^{(0)} Z(\beta^{\left(0\right)},\lambda)}
 \delta_{kj}
- \frac{ \gamma_k^{(0)}(\lambda) \gamma_j^{(0)}(\lambda) }{Z^2(\beta^{\left(0\right)},\lambda)}=Q^{\left(0\right)}_{kj}(\lambda)
\end{multline*}
as $i/m\ri \lambda$. We can then complete the proof of
Proposition \ref{prop1} by applying the dominated convergence theorem
to show the sum converges to the integral (see Section 5.2
of \cite{MPVS} for a similar argument.)

\begin{prop}
\label{prop2}
Under Assumption \ref{A1}
\begin{equation}
\label{2nd_condition}
\lim\limits_{m\ri \infty}\frac{1}{m^2}\Var\left(\sum\limits_{i=2}^{\ell_m}
(\eta_{i}^m)^2\right)=0,
\end{equation}
where the expectation is taken with respect to measure $\P_m^{\left(0\right)}$.
\end{prop}
{\it Proof.}
To simplify  notation we assume in the proof that  $N=1$; 
modifications for the multivariate case are obvious. 
Also, for simplicity of notation, we omit  the upper index in notation for 
$\eta, \zeta$ and $\xi$ variables. So, in the rest of the proof we denote
$\beta=\beta_1,\,  a={\bf a}\in \R,\, \eta_i=\eta_i^m,\, \zeta_i=\zeta^m_{1,i},\,
 \xi_i=\xi_{1,i}^m,\, \bar{\xi}_i=\bar{\xi}^{\, m}_{1,i},\, {\cal F}_{j}=
{\cal F}_{j}^{\left(m\right)}$. 
Besides, we write $\E$ instead of $\E^{\left(0\right)}_m$.

It suffices to show that under Assumption \ref{A1}
\begin{equation}
\label{iii}
\Cov\left(\eta_{i}^2, \eta_{j}^2\right)
\ri 0
\end{equation}
for any pair of sequences $i=i_m$ and $j=j_m$ such that $i\neq j$ and 
 $i/m \ri \lambda',\, j/m\ri \lambda''$ as $m\ri \infty$, 
where $\lambda'$ can coincide with $\lambda''$.
This suffices because the contribution from terms with $i=j$,
divided by $m^2$, is asymptotically negligible since the $\eta_i$
are uniformly bounded.

Recall that  $ \eta_{i}=a(\xi_i-\xbi)/\beta$, where  $\bar{\xi}_{i}=\E\left(\xi_{i}|
{\cal F}_{i-1}\right)$. Therefore,  we need to prove 
that  
\begin{equation}
\label{iv}
\Cov\left((\xi_{i}-\xbi)^2, (\xi_{j}-\xbj)^2\right) \ri 0
\end{equation}
under the same assumptions about the  index sequences.
Assuming for definitness that $i<j$, we have the following identities
\begin{align}
\xi_i^2&=\xi_i, \label{0505a} \\
\E(\xi_i)&=\E\left(\E\left(\xi_i|{\cal F}_{i-1}\right)\right)=\E(\xbi), \label{0505b}\\
\E(\xi_i\xbi)&=\E\left(\xbi\E\left(\xi_i|{\cal F}_{i-1}\right)\right)
=\E(\xbi^2), \label{0505c}\\
\E(f(\xi_i,\xbi, \xbj)\xi_j)&=\E\left(f(\xi_i,\xbi, \xbj)\E\left(\xi_j|{\cal F}_{j-1}\right)\right)=\E(f(\xi_i,\xbi,\xbj)\xbj),
\label{0505d}
\end{align}
where $f(\xi_i,\xbi,\xbj)$ is a polynomial function, e.g., $\xi_i\xbi^2$ etc.
(note that (\ref{0505d}) fails for $i=j$.)
We can write 
$\Cov\left((\xi_{i} - \bar{\xi}_i)^2, (\xi_{j}- \bar{\xi}_j)^2\right)$
as a linear combination of terms of the form
\begin{equation}
\E (\bar{\xi}_i^p \xi_i^{2-p} \bar{\xi}_j^q \xi_j^{2-q})
- \E (\bar{\xi}_i^p \xi_i^{2-p}) 
\E( \bar{\xi}_j^q \xi_j^{2-q})
\label{0505e}
\end{equation}
where $p \in \{0,1,2\}$ and $q \in \{0,1,2\}$.
 As  mentioned before (see display (\ref{xi_conv})), we have that 
\begin{equation}
\bar{\xi}_{r} \toPP b(\lambda):= \frac{\beta_1^{\left(0\right)}\gamma_{1}^{\left(0\right)}(\lambda)}{Z\left(\beta^{\left(0\right)},
 \lambda\right)}
\label{0506a}
\end{equation}
as $r/m\ri \lambda$, 
and also,
 $\E(\xi_{r})\ri b(\lambda)$ as $n\ri \infty$.
Since $\xi_i$ and $\xi_j$ are bounded, we have
$\E(\bar{\xi}_i^2) \to b^2(\lambda')$
and (using  (\ref{0505c}))
$\E(\xi_i \bar{\xi_i}) \to b^2(\lambda')$, while (using (\ref{0505a}))
$\E(\xi^2_i) \to b(\lambda')$, and likewise for $j$. Therefore
\begin{equation}
 \E (\bar{\xi}_i^p \xi_i^{2-p}) \E( \bar{\xi}_j^q \xi_j^{2-q})
\to b^{1 + \min(p,1)}(\lambda') 
 b^{1 + \min(q,1)} (\lambda'').
\label{0505f}
\end{equation}
But
 using (\ref{0505d}), (\ref{0506a})  
and Proposition \ref{dct_prob}
we also have
\begin{equation}
 \E(\bar{\xi}_i^p \xi_i^{1-p} \bar{\xi}_j \xi_j ) =
\E(\bar{\xi}_i^p \xi_i^{1-p} \bar{\xi}_j^2) 
 \to b^{1+ \min(p,1)}(\lambda') b^2(\lambda''),
\label{0506b}
\end{equation}
and
 using (\ref{0505a}) for $j$,  (\ref{0505d}), (\ref{0506a})
and Proposition \ref{dct_prob}
we also have
\begin{equation}
 \E(\bar{\xi}_i^p \xi_i^{1-p}  \xi^2_j )
 = \E(\bar{\xi}_i^p \xi_i^{1-p}  \bar{\xi}_j )
 \to b^{1+ \min(p,1)}(\lambda') b(\lambda'').
\label{0506c}
\end{equation}
Combining (\ref{0506b}) and (\ref{0506c}) shows that
$
\E(\bar{\xi}_i^p \xi_i^{1-p} \bar{\xi}_j^q \xi_j^{1-q})
$ converges to the same limit as the expression in (\ref{0505f}). 
Hence, each expression of the form in (\ref{0505e}) tends to zero,
and we have established (\ref{iv}). Hence,
Proposition \ref{prop2} is proved.

\section{Numerical example}
\label{example}

In this section we give a numerical example 
demonstrating 
that MLE is effective in distinguishing  between CSA's which might 
generate  similar patterns.
 
\begin{figure}[htbp]
\centering
\begin{tabular}{cc}
    \includegraphics[width=2.5in, height=2.5in, angle=270]{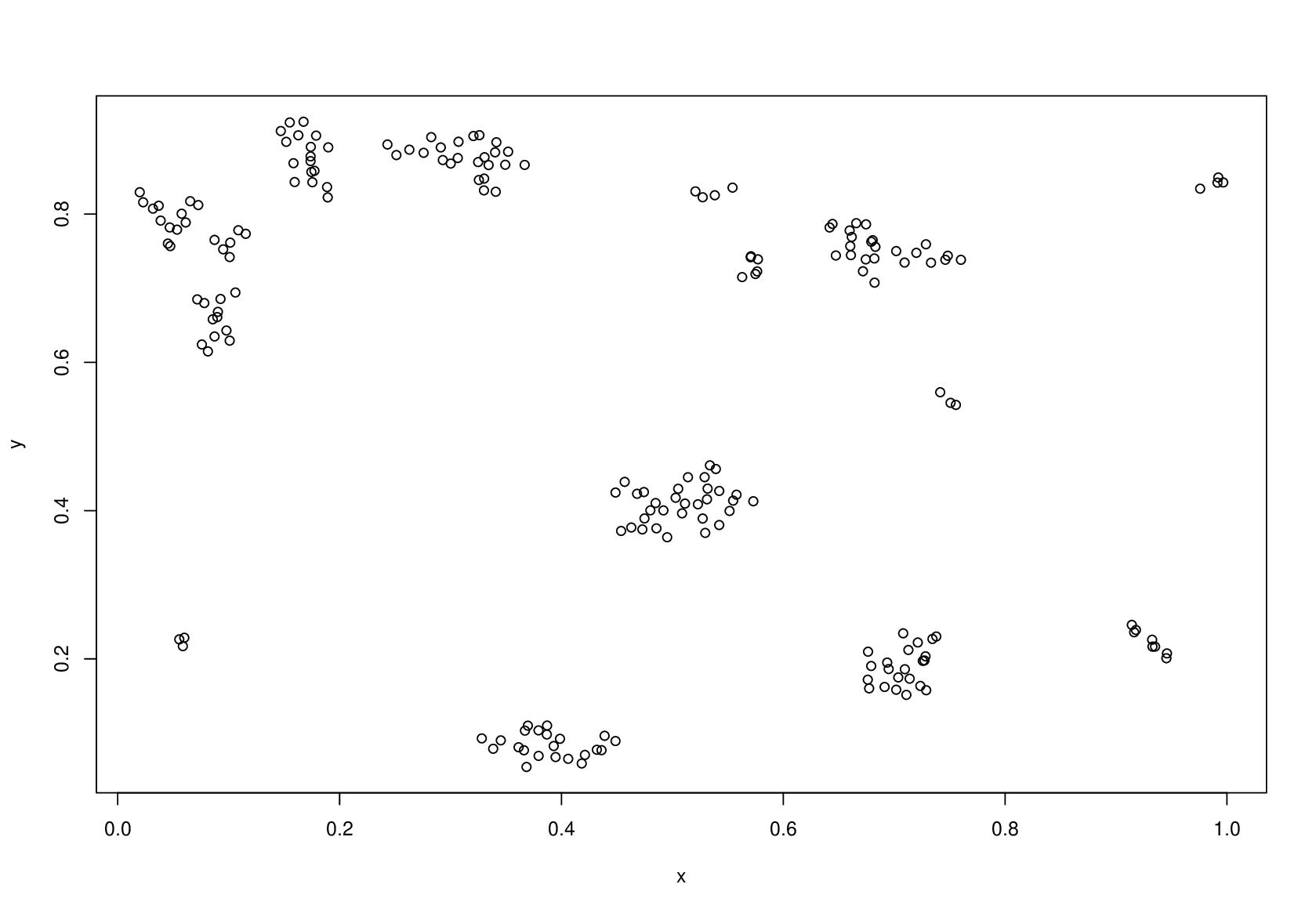} &
    \includegraphics[width=2.5in, height=2.5in, angle=270]{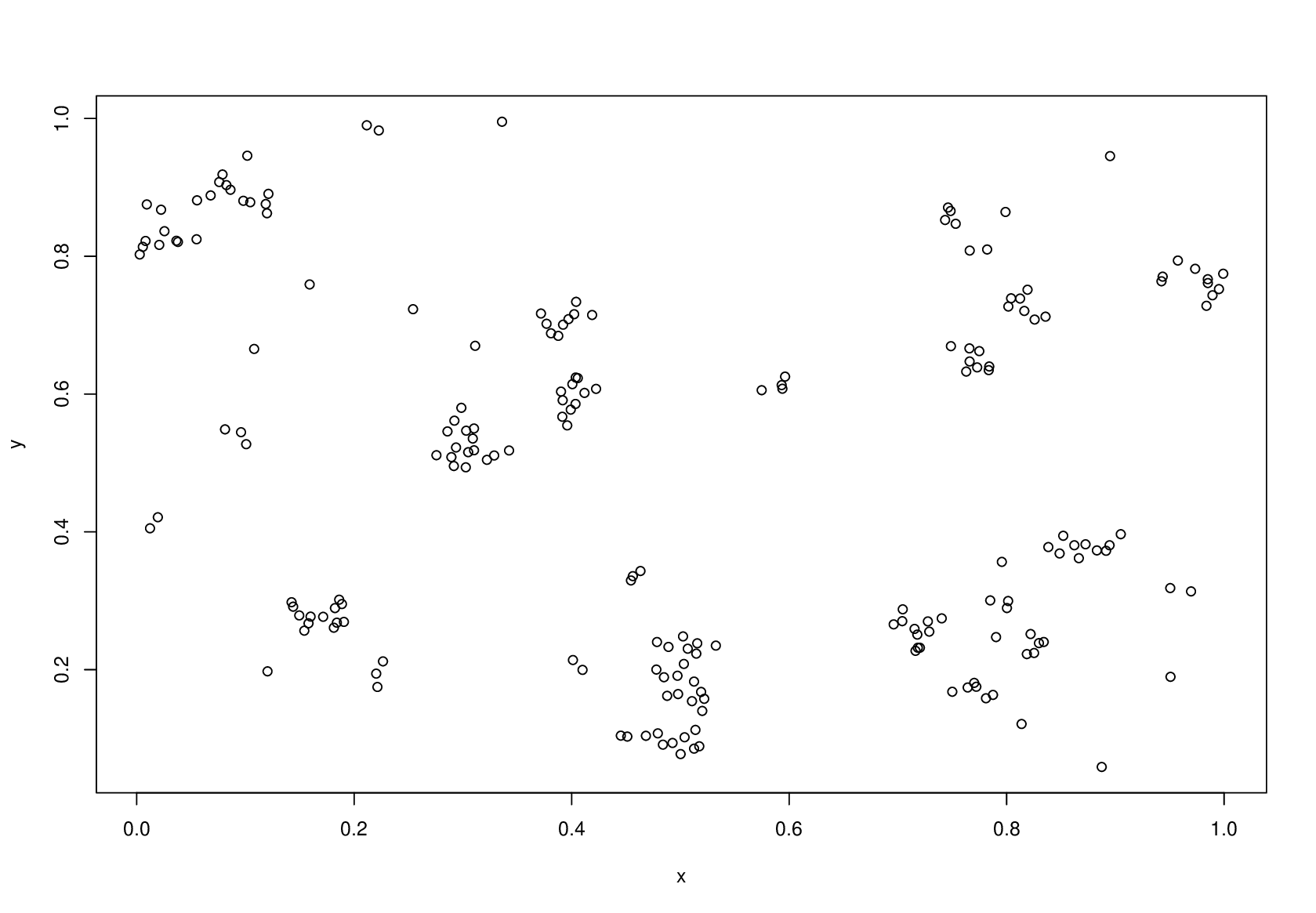} \\
\end{tabular}
\caption{{\footnotesize $\ell=200$. Left: increasing rates,
$(t_0,t_1,t_2)=(16, 93, 91)$.
Right: flat rates,  $(t_0,t_1,t_2)=(43, 100, 57)$.
}}
\label{fig21}
\end{figure}

\begin{figure}[htbp]
\centering
\begin{tabular}{cc}
    \includegraphics[width=2.5in, height=2.5in, angle=270]{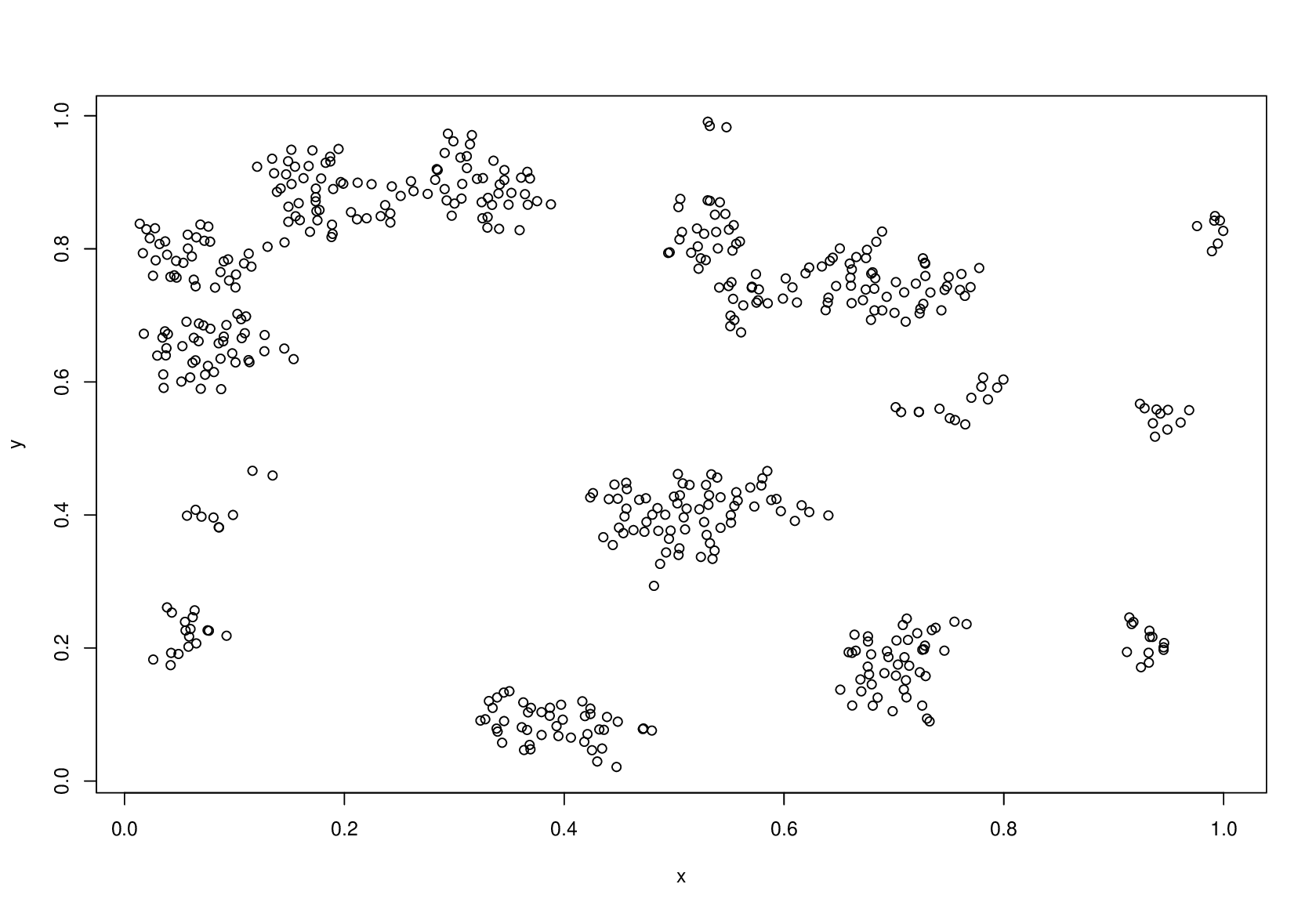} &
    \includegraphics[width=2.5in, height=2.5in, angle=270]{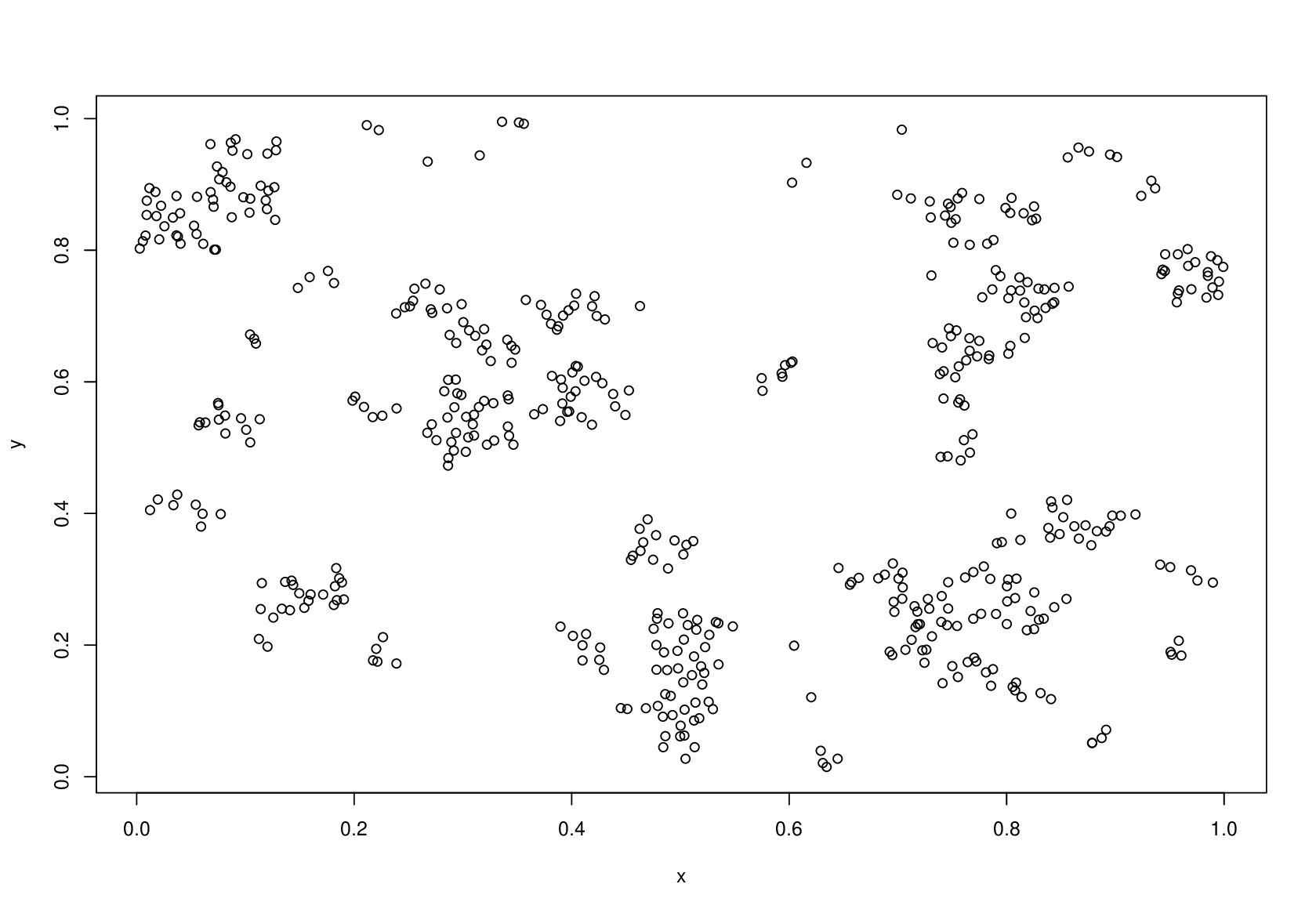} \\
\end{tabular}
\caption{{\footnotesize $\ell=500$. Left: increasing rates,
$(t_0,t_1,t_2)=(25, 233, 242)$.
Right: flat rates, $(t_0,t_1,t_2)=(62, 272, 166)$.}}
\label{fig22}
\end{figure}

\begin{figure}[htbp]
\centering
\begin{tabular}{cc}
    \includegraphics[width=2.5in, height=2.5in, angle=270]{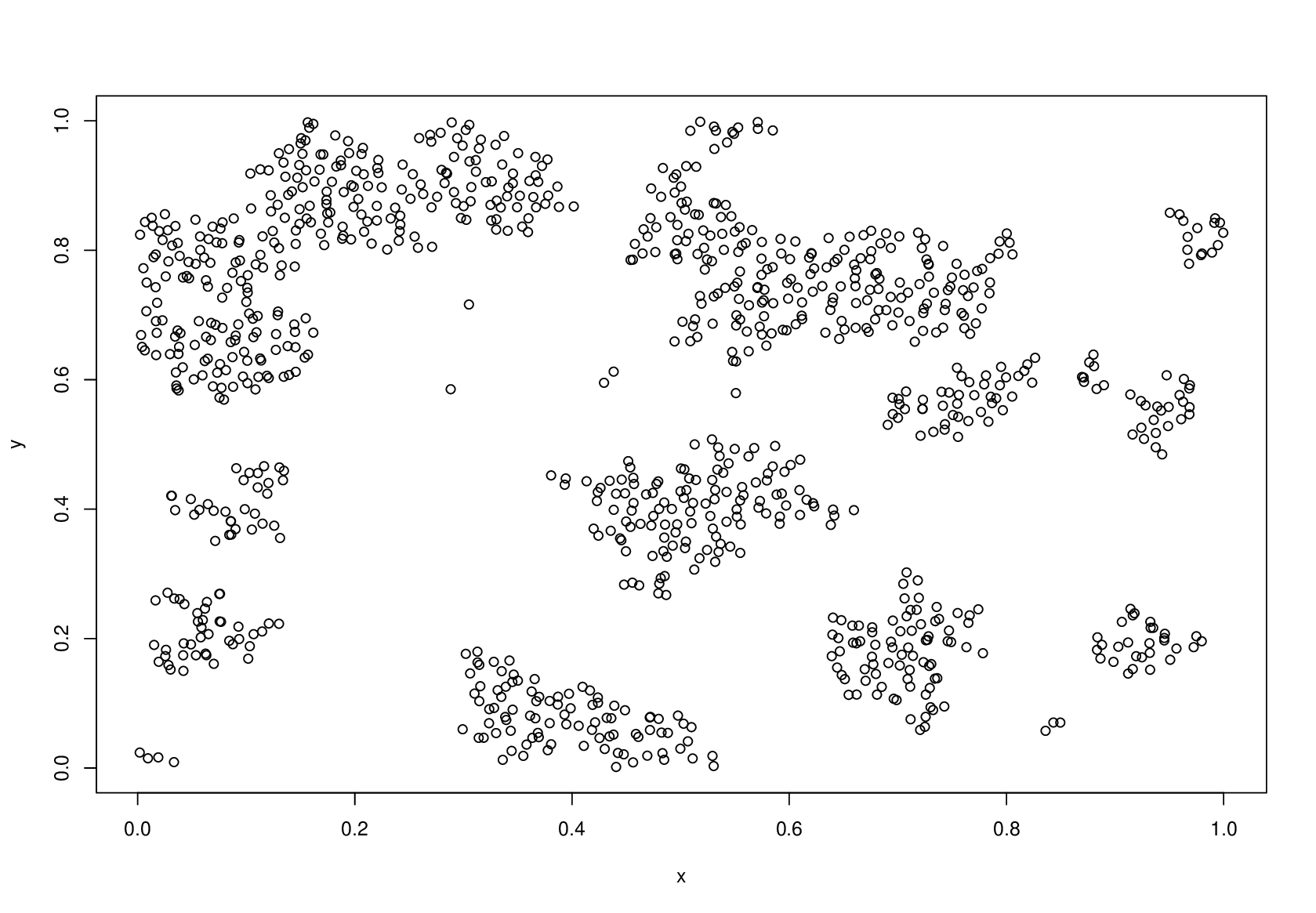} &
    \includegraphics[width=2.5in, height=2.5in, angle=270]{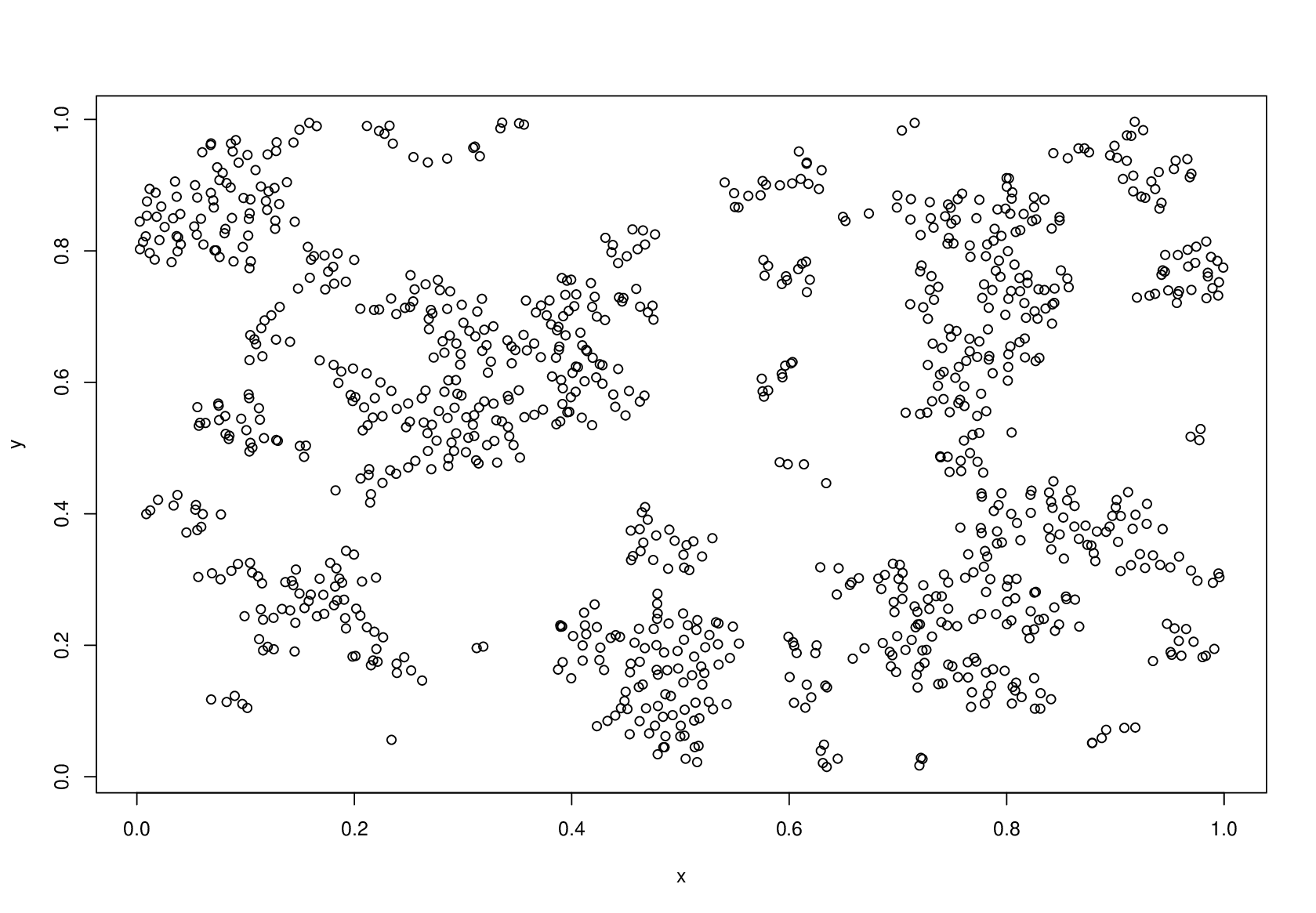} \\
\end{tabular}
\caption{{\footnotesize $\ell=1000$. 
Left: increasing rates,  $(t_0,t_1,t_2)=(34, 434, 532)$.
Right: flat rates, 
$(t_0,t_1,t_2)=(84, 552, 364)$.
}}
\label{fig23}
\end{figure}

\begin{figure}[htbp]
\centering
\begin{tabular}{cc}
    \includegraphics[width=2.5in, height=2.5in, angle=270]{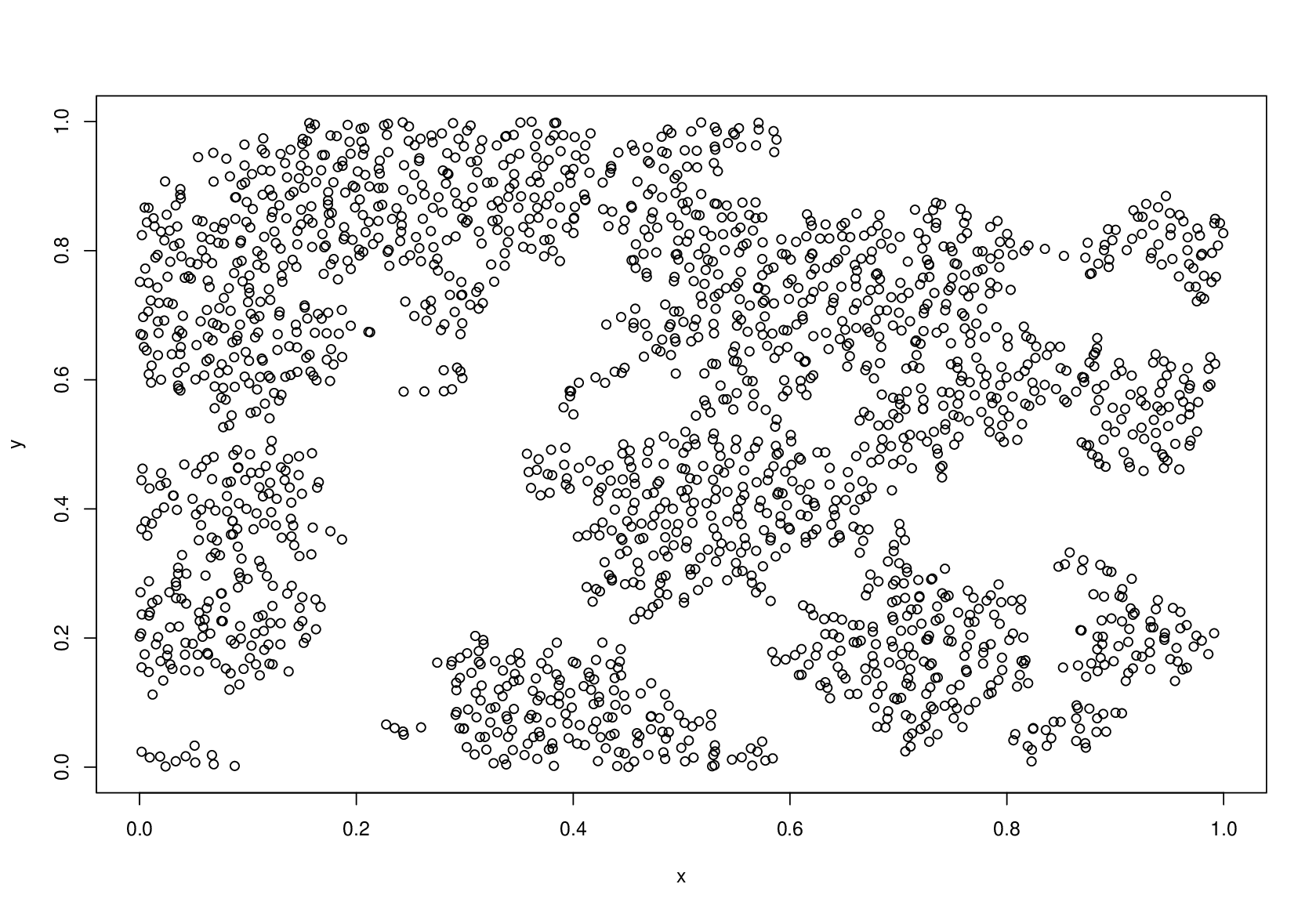} &
    \includegraphics[width=2.5in, height=2.5in, angle=270]{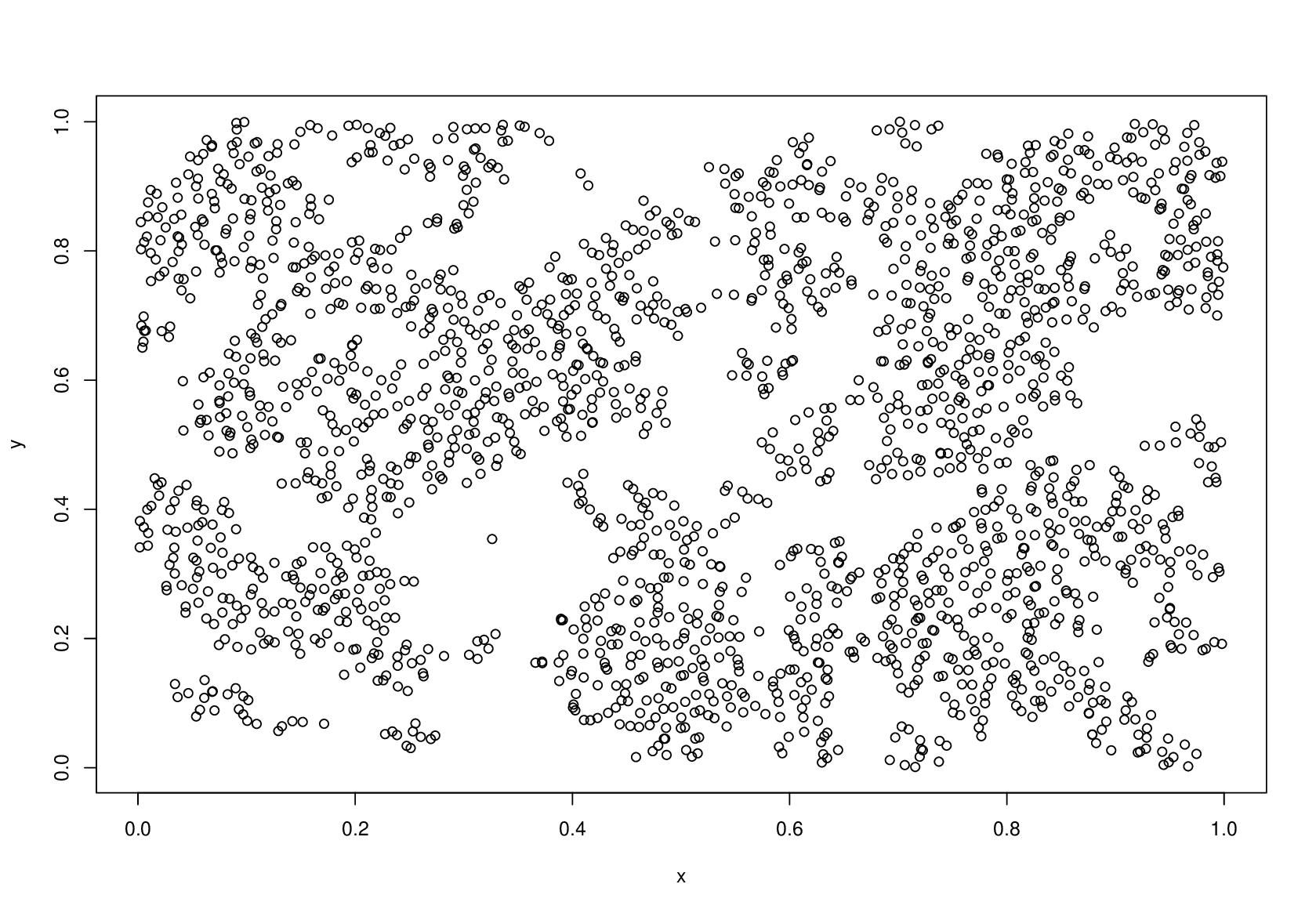} \\
\end{tabular}
\caption{{\footnotesize $\ell=2000$. 
Left: increasing rates,
$(t_0,t_1,t_2)=(43, 825, 1132)$.
Right: flat rates, $(t_0,t_1,t_2)=(95, 1048, 857)$.
}}
\label{fig24}
\end{figure}

\begin{figure}[htbp]
\centering
\begin{tabular}{cc}
    \includegraphics[width=2.5in, height=2.5in, angle=270]{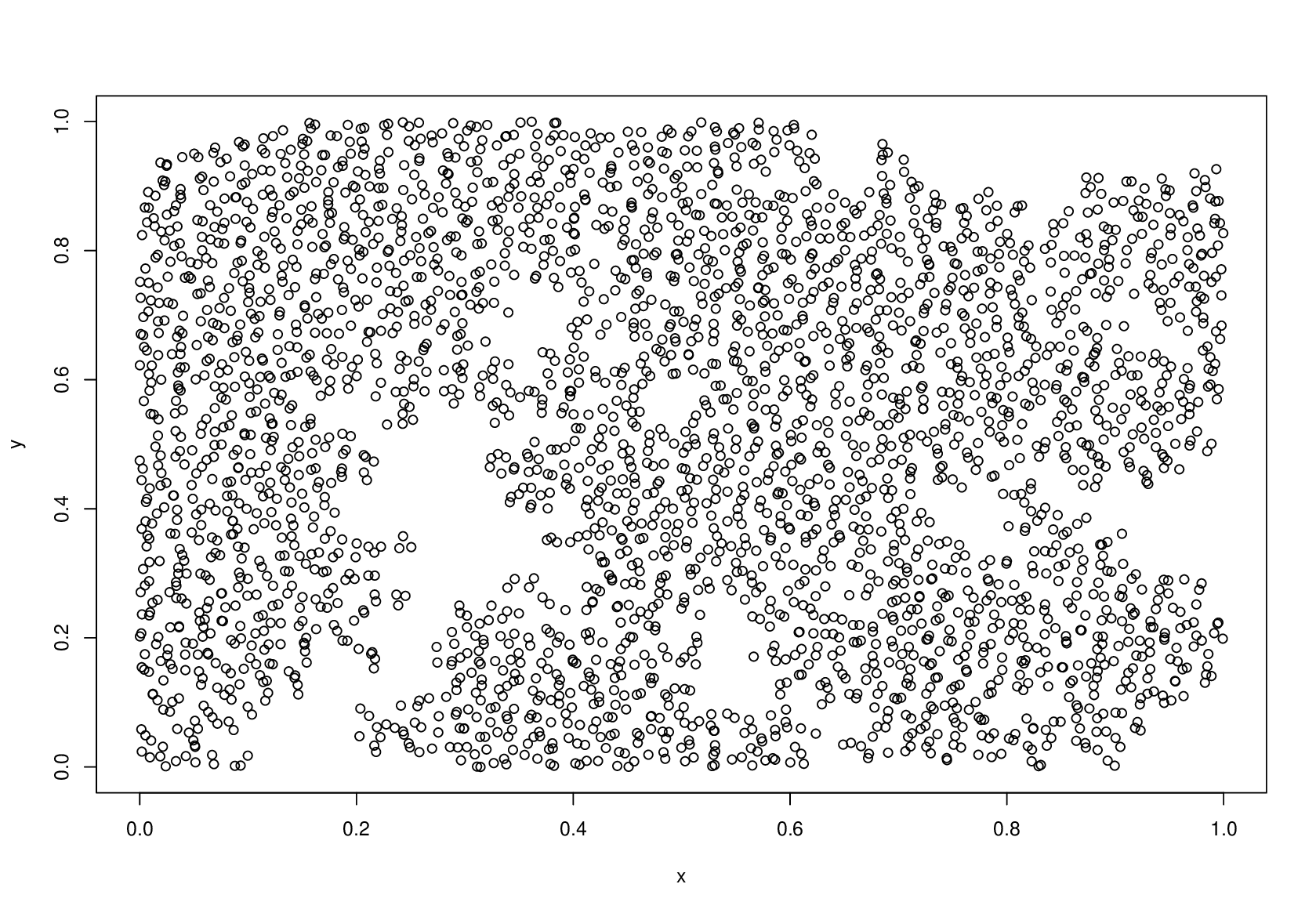} &
    \includegraphics[width=2.5in, height=2.5in, angle=270]{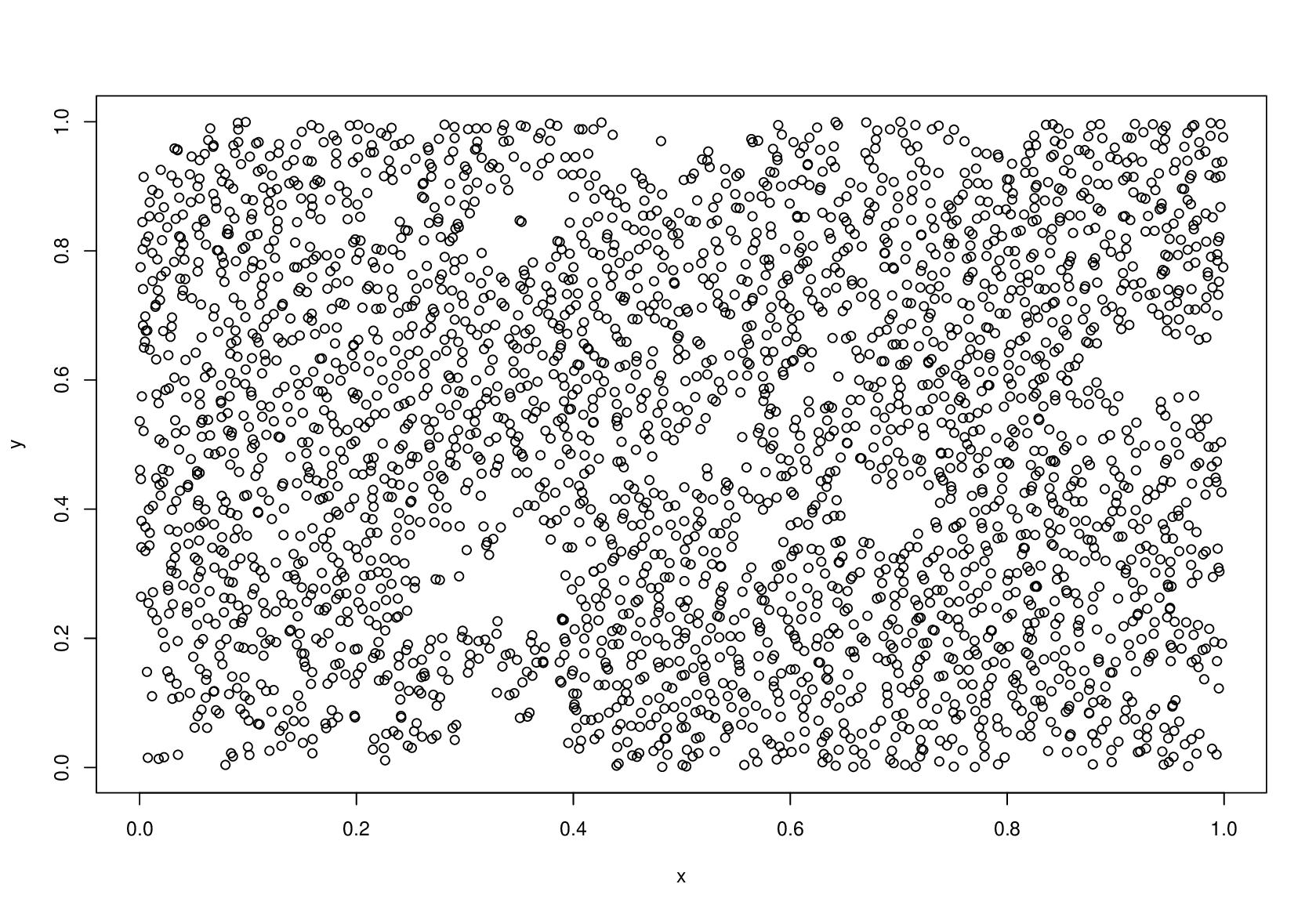} \\
\end{tabular}
\caption{{\footnotesize $\ell=3000$
Left: increasing rates, $(t_0,t_1,t_2)=(47, 1190, 1763)$.
Right:  flat rates,  $(t_0,t_1,t_2)=(106, 1473, 1421)$.
}}
\label{fig25}
\end{figure}

\begin{figure}[htbp]
\centering
\begin{tabular}{cc}
    \includegraphics[width=2.6in, height=2.6in, angle=270]{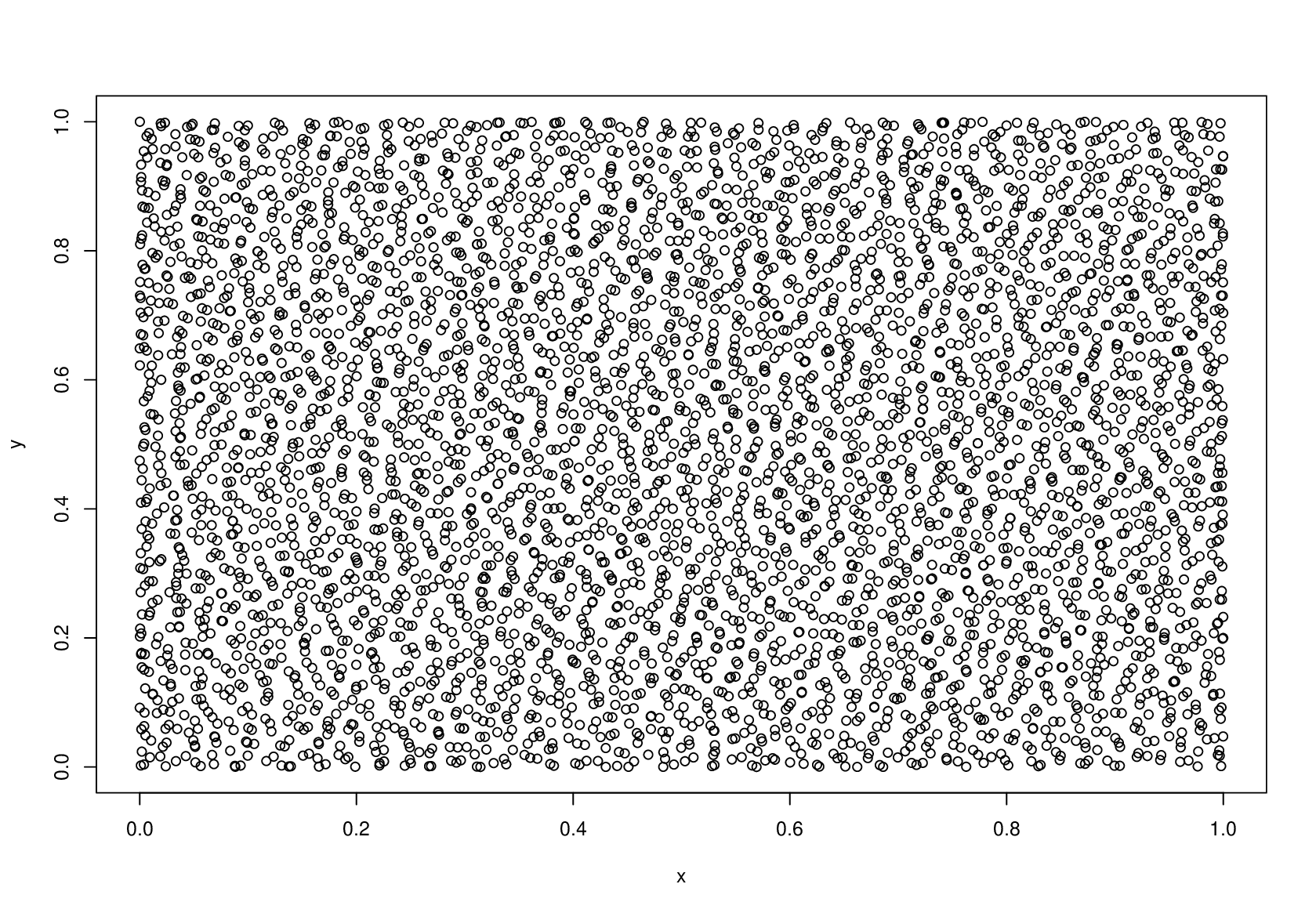} &
    \includegraphics[width=2.6in, height=2.6in, angle=270]{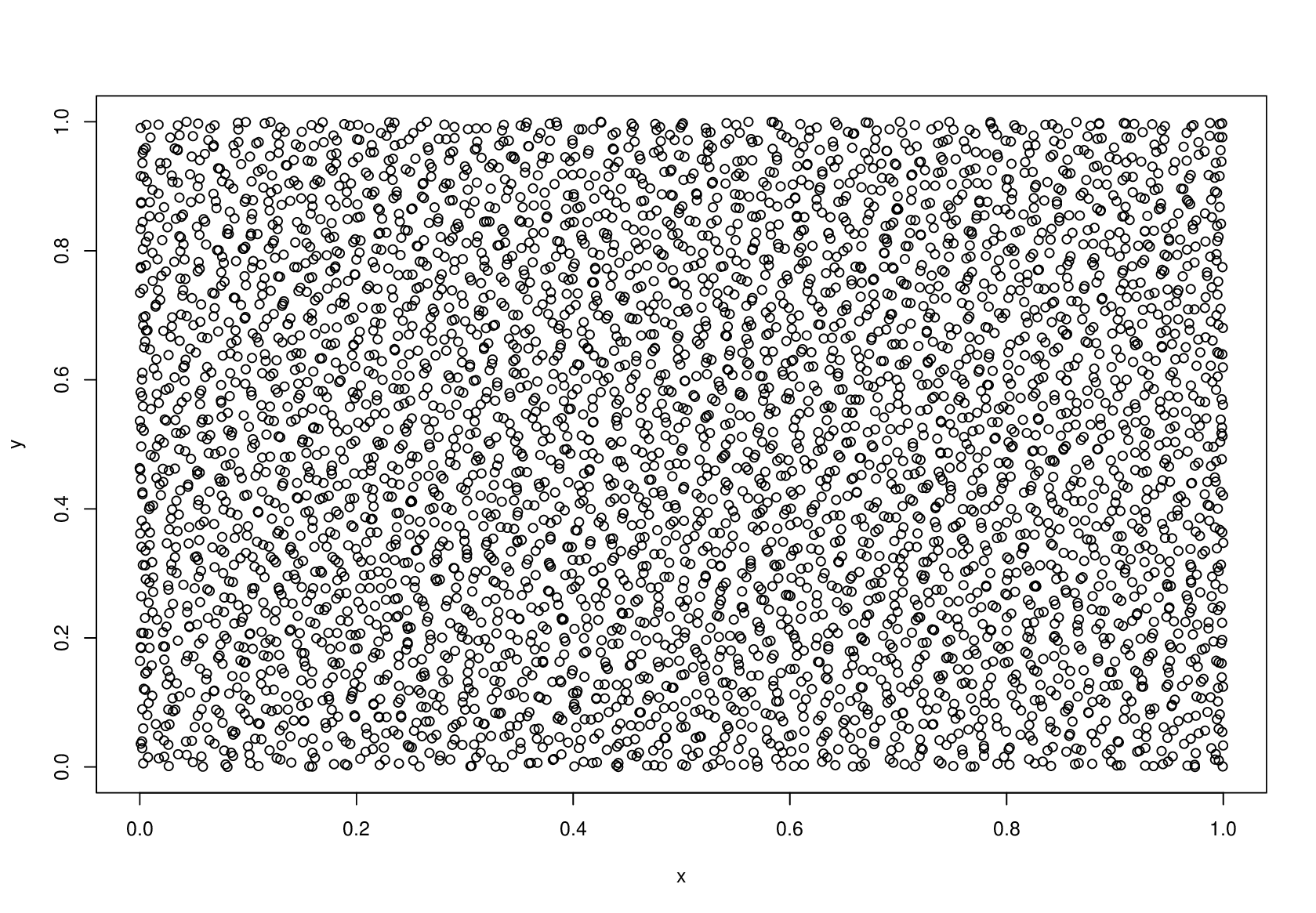} \\
\end{tabular}
\caption{{\footnotesize 
Left: increasing rates,$\ell=4407$, $(t_0,t_1,t_2)=(48, 1426, 2933)$.
Right:  flat rates, $\ell=4416$,   $(t_0,t_1,t_2)=(108, 1688, 2620)$.
}}
\label{fig26}
\end{figure}

In  \cite{MPVS} we briefly discussed difference between  clustering effects produced by 
CSA determined by  a set of increasing  parameters $\beta$ (the so-called Aarhenius rates, 
\cite{Evans}), 
and determined by  a set of flat rates (the so-called Eden rates, \cite{Evans}).
 As before, we consider two single realizations of CSA. 
  Six successive  images for each of realization shown in Figures \ref{fig21}-\ref{fig26}. 
The interaction radius is $R=0.02$ in both cases.
The left images have been generated by CSA  with $\beta$-parameters $\beta_0=1, 
\beta_1=300, \beta_2=500, \beta_k=0, \, k\geq 3$. The right images
have been generated by CSA with $\beta$-parameters
 $\beta_0=1, \beta_1=\beta_2=100, \beta_k=0\, k\geq 3$. 
The first five pairs of images 
 with first $\ell=200, 500, 1000, 2000$ and $\ell=3000$ points respectively are shown in Figures \ref{fig21}-\ref{fig25}. The last pair of images shows the realisations at jamming, i.e.,
when there is no space left to accommodate a point. The left image contains $\ell=4407$ and 
the right image contains $\ell=4416$ points. 
Can one tell apart these two sets of parameters given the series of images provided?

The images with $200$ points look similar and it seems plausible that they 
have been generated by the same CSA. 
In both cases
  new points  tend  to appear in the vicinity of existing points because of the choice 
of the parameters. Though clusters formed by a single point are noticeable on the right image
and clusters seem to be more dense on the left one.

The  pair of subsequent  images containing  $500$ points is shown
 in Figure \ref{fig22}. 
It is noticeable at both images there are almost no  new  clusters; the existing 
clusters keep growing and eventually start coalescing. 
Besides, 
it is slightly  visible that the right pattern is more dispersed than the left one.
All these effects  are becoming  more  visible  for the  pair of images  showing 
further evolution and containing  $1000$ and $2000$ points. These images  
are shown in Figures \ref{fig23} and \ref{fig24}.

 The effects that have been just described 
 are rather straightforward analogues of the 
 phenomenon of
``competition between the birth, growth and coalescence''
 (\cite{Evans}, p.1307),
 which is well known for lattice CSA models.

 Though the main basic feature of both
 series of  images, namely, {\it clustering},
is common to both choices of the parameters, 
 the clustering effect is  more visible
 in the images produced 
by  the model with an increasing set of non-zero parameters  
(the sequence of left images).
 The clusters are more saturated in the left images, i.e. clustering 
is stronger. It seems that the right realisation spreads faster in comparison to the left one.
This is called {\it mild} clustering; the distribution  of points inside a cluster is more or less regular, since 
a new  point distribution is uniform  conditioned on being adsorbed in the vicinity of existing points.

The difference between the strong and the  mild clustering (corresponding  to  increasing
and flat sets of non-zero parameters respectively) observed in Figures \ref{fig22}-\ref{fig24},  
vanishes at the later stages of evolution,
when it approaches jamming.
It is quite difficult  to distinguish by visual inspection the two sets of parameters given 
the pair of images shown in Figure  \ref{fig25}. 
Note that both of these images are close to the corresponding
 jamming images shown in Figure \ref{fig26}.
One might argue that these two realisations have been produced by the same model
and the  differences between them (observed at some intermideate images) 
can be attributed to variability of the samples. 
Numerical results given in Tables \ref{T1} and \ref{T2} 
show that MLE is an effective tool for parameter estimation. The tables
contain MLE's for both sets of parameters along with corresponding  approximate confidence bounds
(any computed value is rounded  to its  nearest integer).
The $95\%$ confidence bounds are computed  by {\it  formally assuming} normality of 
 $\widehat{\beta}$.
The variances of the estimates are approximated, as usual,  by the corresponding diagonal elements
of  the matrix inverse to the observed information matrix. The latter turned out to be non-degenerate
for all observed images.
The variances of the estimates  decrease as the number of observed points increase.
As a result, the confidence intervals become narrower. 
The tendency breaks down only for  the rightmost entry 
of  the bottom line in Table \ref{T1}. Perhaps this can be 
explained by the lack of accuracy of the computations (see the discussion 
of computational issues in \cite{MPVS}). 
The observed reduction of variances  is intuitively expected, 
although  the normality 
assumption in the unit volume cannot be based on our asymptotic results. 
This is in contrast to the limiting situation where the effect 
is clearly implied by  the integral representation (\ref{integ}) for the information matrix.
The representation implies that the variance of the estimate $\widehat{\beta}_i$ converges,
as $m\ri \infty$ and $\ell_m/m\ri \mu$, to
$$\frac{1}{\int_{0}^{\mu} g_{ii}(\lambda) d\lambda},$$ 
where $g_{ii}(\lambda)\geq 0$ is the $i-$th eigenvalue of  matrix $ Q^{\left(0\right)}(\lambda)$
in the representaion (\ref{integ}). The preceeding display justifies 
 "reduction of   variances" effect, if the density of points, i.e. $\mu$, increases. 
The lower bound for the variance of the estimate $\widehat{\beta}_i$ 
is given by the same formula with $\mu=\theta_{\infty}$,
where $\theta_{\infty}$ is the jamming density (\cite{MPVS}).

Under certain assumptions normality  of $\widehat{\beta}$ 
in a {\it fixed finite} volume can {\it possibly} be  advocated as follows.
Consider, for definiteness, the model in the unit volume
and let the interaction radius be sufficiently small. This is the case in the simulated 
examples.  If the interaction radius 
is sufficiently small, then the jamming density is high.
In other words, a sufficiently large number of points  can be accommodated.
 It was shown in Section \ref{CLT} that 
the score function is a martingale sum  containing  $\ell$ terms, where 
$\ell$ is the number of observed points.
Therefore, one might expect that if $\ell$ is sufficiently 
large (e.g., thousands), then  the normal approximation starts working.

Finally, it should be noted that  MLEs effectively  
capture
the correct magnitude of the parameters 
and this is why  two considered 
sets of parameters in the example (producing sometimes quite similar images) 
can be effectively  distinguished. 
For the sake of completeness, consider also  the left image in Figure \ref{fig1}. 
It has been  generated by CSA with the interaction radius 
$0.01$ and $\beta-$parameters $\beta_0=1.0, \beta_1=1000.0, \beta_2=10000.0, \beta_k=0.0, \, k\geq 3.$
The image contain $\ell=1000$ points, $t-$statistics are $t_0=23, t_1=149, t_2=828$. 
The MLE estimates for $\beta_1$ and $\beta_2$ are $1105.0$ and $10510.0$ respectively.

\begin{table}[htb]
\caption{{\small MLE's  for the left images in Figures \ref{fig21}-\ref{fig26}}}
\footnotesize
\begin{center}
\begin{tabular}{|c|c|c|c|c|c|}
\hline
  $\ell=200$ &  $\ell=500$  & $\ell=1000$  & $\ell=2000$ & 
 $\ell=3000$ & $\ell=4407$\\
\hline
 $\widehat{\beta}_1=401$ & $\widehat{\beta}_1=377$ & 
$\widehat{\beta}_1=334$ & $\widehat{\beta}_1=320$ & $\widehat{\beta}_1=318$ & 
$\widehat{\beta}_1=323$  \\
 (176, 626) & (214, 540) & (213,  455) &
(218, 422) & (223, 413) & (226,  420)
  \\ 
\hline
 $\widehat{\beta}_2=695$ & $\widehat{\beta}_2=594$ & 
$\widehat{\beta}_2=566$ & $\widehat{\beta}_2=546$ & 
$\widehat{\beta}_2=521$ &$\widehat{\beta}_2=533$ \\
 (298, 1091) & (335, 853) & (360, 772) &
(371, 721) & (364, 678) & (373,  693)  \\ 
\hline
\end{tabular}
\end{center}
\label{T1}
\end{table}

\begin{table}[htb]
\caption{{\small MLE's for the right images in Figures \ref{fig21}-\ref{fig26}}}
\footnotesize
\begin{center}
\begin{tabular}{|c|c|c|c|c|c|}
\hline
  $\ell=200$ &  $\ell=500$  & $\ell=1000$  & $\ell=2000$ & 
 $\ell=3000$ & $\ell=4416$ \\
\hline
 $\widehat{\beta}_1=89$ & $\widehat{\beta}_1=98$ & 
$\widehat{\beta}_1=96$ & $\widehat{\beta}_1=104$ & $\widehat{\beta}_1=101$ 
& $\widehat{\beta}_1=98$\\
(55, 123) & (69, 127) & (73, 119) &
(81, 127) & (80, 122) & (78, 118)\\ 
\hline
$\widehat{\beta}_1=106$ & $\widehat{\beta}_1=97$ & 
$\widehat{\beta}_2=88$ & $\widehat{\beta}_2=100$ & $\widehat{\beta}_2=99$ &
 $\widehat{\beta}_2=98$\\
 (61, 151) & (67, 127) & (66, 110) &
(78, 122) & (78, 120) & (78,  118) \\ 
\hline
\end{tabular}
\end{center}
\label{T2}
\end{table}

\newpage 

\section*{Appendix. On positive definiteness of the limit information matrix}
\label{matr}
It is easy to see from equation (\ref{integ}) that positive definiteness of matrix 
$ Q^{\left(0\right)}(\lambda)=Q\left(\beta^{\left(0\right)}, \lambda\right)$
for any fixed $\lambda\in (0,\theta_\infty)$ 
implies  positive definiteness of the limit matrix $J^{\left(0\right)}(\mu)$.
Positive definiteness of matrix $ Q^{\left(0\right)}(\lambda)$
was shown in Lemma 5.2 in \cite{MPVS}. 
Here we give another proof by studying  the  matrix structure in more detail.

It can be seen from equation (\ref{mat}) that  the matrix principal minor  formed by 
the  intersection of the first $k$ rows and $k$ columns is
\[
D_{N,k}(\beta^{\left(0\right)},\lambda)=
 \left( \begin{array}{ccc}
\frac{\gamma_1^{\left(0\right)}(\lambda)
\left(Z(\beta^{\left(0\right)},\lambda)-
\gamma_1^{\left(0\right)}(\lambda)\beta^{\left(0\right)}_1\right)}
{\beta^{\left(0\right)}_1 Z^2(\beta^{\left(0\right)},\lambda)} & 
\ldots &  -\frac{\gamma_1^{\left(0\right)}(\lambda) \gamma_k^{\left(0\right)}(\lambda)}
{Z^2(\beta^{\left(0\right)},\lambda)}\\
\vdots & \ldots  & \vdots  \\
 -\frac{\gamma_1^{\left(0\right)}(\lambda) \gamma_k^{\left(0\right)}(\lambda)}
{Z^2(\beta^{\left(0\right)},\lambda)} & \ldots &
\frac{\gamma_k^{\left(0\right)}(\lambda)
\left(Z(\beta^{\left(0\right)},\lambda)-
\gamma_k^{\left(0\right)}(\lambda)\beta^{\left(0\right)}_k\right)}
{\beta^{\left(0\right)}_k Z^2(\beta^{\left(0\right)},\lambda)}
\end{array} \right).
\] 
It is easy to see that  determinant of $D_{N,k}(\beta^{\left(0\right)},\lambda)$ is  
$$
\left|D_{N,k}(\beta^{\left(0\right)},\lambda)\right|=
\frac{(-1)^k}
{Z^{2k}(\beta^{\left(0\right)},\lambda)}
\prod\limits_{i=1}^k\frac{\gamma_i^{\left(0\right)}(\lambda)}
{\beta_i^{\left(0\right)}} \left|A_k-Z(\beta^{\left(0\right)},\lambda)E_k\right|,$$
where $|A_k-Z(\beta^{\left(0\right)},\lambda)E_k|$ is  determinant of  matrix 
$A_k-Z(\beta^{\left(0\right)},\lambda)E_k$, where, in turn,   matrix $A_k$ is defined as follows
\begin{equation}
\label{matrix}
A_{k}=\left(\beta^{\left(0\right)}_1,\ldots,\beta^{\left(0\right)}_k\right)
(\gamma_1^{\left(0\right)}(\lambda),\ldots,\gamma_k^{\left(0\right)}(\lambda))^T,
\end{equation}
and 
$E_k$ is the $k\times k$ unit matrix.
By definition, $|A_k-Z(\beta^{\left(0\right)},\lambda)E_k|$ is the characteristic 
polynomial of $A_k$ evaluated at point $Z(\beta^{\left(0\right)},\lambda)$.
It can be shown (we omit the proof) 
 that if  $a,b\in \Cc^n$ are non-zero complex  vectors, such that $a^Tb\neq 0$, 
then a  quadratic matrix  $M=ab^T$  
has the only non-zero eigenvalue $a^Tb$ of multiplicity $1$, 
$0$ is the other  matrix eigenvalue of multiplicity $n-1$ and
the matrix characteristic polynomial is 
$$|M-uE|=(-1)^nu^{n-1}\left(u-a^Tb\right),\,\,\ u\in \Cc^n.$$
Hence,   
\begin{align*}
|A_k-Z(\beta^{\left(0\right)},\lambda)E_k|
&= (-1)^k Z^{k-1}(\beta^{\left(0\right)},\lambda)
\left(\gamma^{\left(0\right)}_0(\lambda)+
\sum\limits_{i=k+1}^N\beta^{\left(0\right)}_i\gamma^{\left(0\right)}_i(\lambda)\right)
\end{align*}
and  
$$
\left|D_{N,k}(\beta^{\left(0\right)},\lambda)\right|=
\frac{\left(\gamma^{\left(0\right)}_0(\lambda)+
\sum_{i=k+1}^N\beta^{\left(0\right)}_i
\gamma^{\left(0\right)}_i(\lambda)\right)}{Z^{k+1}(\beta^{\left(0\right)},\lambda)}
\prod\limits_{i=1}^k\frac{\gamma_i^{\left(0\right)}(\lambda)}
{\beta_i^{\left(0\right)}}.
$$ 
The right  side of the preceding display is  positive because 
the functions $\gamma_i,\, i=1,\ldots,N$ are positive. 
Thus any principal minor of matrix (\ref{mat}) is positive
and  by Sylvester criterion this  matrix is positive definite.

\end{document}